\documentclass[reqno]{amsart}
\usepackage{amsfonts,amssymb,amsmath,amscd}
\usepackage{graphicx}
\usepackage{psfrag}
\usepackage{xspace}
\usepackage{colordvi}

\date{25 October 2003 (Journal submission 15 December 2000)}

\title[Interlace polynomial]{
The Interlace Polynomial of a Graph }

\author[Arratia]{Richard Arratia}
\address[Richard Arratia]{Univ.\ of Southern California \\
Department of Mathematics \\
Los Angeles CA 90089, USA}
\email{rarratia@math.usc.edu}

\author[Bollob\'{a}s]{B\'{e}la Bollob\'{a}s$^\dag$}
\address[B\'{e}la Bollob\'{a}s]{Univ.\ of Memphis \\
Department of Mathematical Sciences \\
Memphis TN 38152, USA,
and
Trinity College, Cambridge CB2 1TQ, U.K.
}
\email{bollobas@msci.memphis.edu}

\author[Sorkin]{Gregory B. Sorkin}
\address[Gregory B. Sorkin]{IBM T.J.\ Watson Research Center \\
Department of Mathematical Sciences \\
Yorktown Heights NY 10598, USA}
\email{sorkin@watson.ibm.com}

\thanks{$^\dag$
Bollob\'{a}s was supported by NSF grant DMS-9971788.
}


\begin{document}
\bibliographystyle{amsalpha}

\newtheorem{theorem}{Theorem}
\newtheorem{conjecture}[theorem]{Conjecture}
\newtheorem{corollary}[theorem]{Corollary}
\newtheorem{definition}[theorem]{Definition}
\newtheorem{lemma}[theorem]{Lemma}
\newtheorem{proposition}[theorem]{Proposition}
\newtheorem{remark}[theorem]{Remark}
\newtheorem{example}[theorem]{Example}

\newcommand{\csub}[1]{{\makebox[0in]{$\substack{#1}$}}}
\newcommand{\Gstar}{G^\star}
\newcommand{\Gab}{G^{ab}}
\newcommand{\Hab}{H^{ab}}
\newcommand{\uu}{a}
\newcommand{\vv}{b}
\newcommand{\uv}{ab}
\newcommand{\vu}{ba}
\newcommand{\by}{\times}
\newcommand{\qq}{\tilde{q}}
\newcommand{\ilace}{H}
\newcommand{\HH}{H}
\newcommand{\Hset}{\mathcal{H}}
\newcommand{\Pset}{\mathcal{P}}
\newcommand{\Cset}{\mathcal{C}}
\newcommand{\Dset}{\mathcal{D}}
\newcommand{\Dstar}{D^{\star}}
\newcommand{\Cstar}{C^{\star}}
\newcommand{\Dstarp}{{D^{\star}}'}
\newcommand{\Cstarp}{{C^{\star}}'}
\newcommand{\qh}{q_\HH}
\newcommand{\qhx}{q_\HH(x)}
\newcommand{\R}{\mathbb{R}}
\def\tito{2-in, 2-out }
\newcommand{\ZZ}{\mathbb Z}
\def\cal{\mathcal}
\newcommand{\ckt}{C}
\newcommand{\order}[1]{|{#1}|}
\newcommand{\numeuler}{r_1}
\newcommand{\size}[1]{e({#1})}
\newcommand{\ie}{\textit{i.e.}}
\newcommand{\eg}{\textit{e.g.}}
\newcommand{\restricted}[1]{\: |_{#1}}
\newcommand{\eqdef}{\stackrel{\mbox{\scriptsize\rm def}}{=}}
\newcommand{\ind}{\hspace*{1cm}}
\newcommand{\tmbf}[1]{\textbf{\boldmath{{#1}}}}
\newcommand{\isom}{\cong}
\newcommand{\vleq}{\preceq}
\newcommand{\vgeq}{\preceq}

\begin{abstract}
Motivated by circle graphs, and the enumeration of Euler circuits, we define a one-variable ``interlace polynomial'' for any graph.  The polynomial satisfies a beautiful and unexpected reduction relation, quite different from the cut and fuse reduction characterizing the Tutte polynomial.

It emerges that the interlace graph polynomial may be viewed as a special case of the Martin polynomial of an isotropic system, which underlies its connections with the circuit partition polynomial and the Kauffman brackets of a link diagram.  The graph polynomial, in addition to being perhaps more broadly accessible than the Martin polynomial for isotropic systems, also has a two-variable generalization that is unknown for the Martin polynomial.  We consider extremal properties of the interlace polynomial, its values for various special graphs, and evaluations which relate to basic graph properties such as the component and independence numbers.
\end{abstract}

\maketitle
\tableofcontents

\section{Introduction}
\subsection{Motivation}
This work was originally motivated by a problem relating to
DNA sequencing by hybridization~\cite{ABCS00}.
At the mathematical heart of the problem was to count the
number of \tito digraphs (Eulerian, directed graphs in which
each vertex has in-degree two and out-degree two) having a
given number of Euler circuits,
formulas for which are given in~\cite{ABCS00}.

This may be thought of as a sort of inverse of a standard problem,
counting the Euler circuits in a directed graph.
That can be done (in polynomial time) with a combination of
the matrix-tree theorem~\cite{KirchoffMT,Tutte48}
(see also~\cite[p.~58]{BolMGT}) and
the so-called BEST theorem~\cite{ST41,BE51} (see also~\cite[p.~18]{BolMGT}).
However, this approach does not give the structural information that
was needed to solve the inverse problem.

The successful approach led to the graph polynomial discussed here.
In addition to counting Euler circuits in a \tito digraph,
the interlace polynomial
(like the Martin polynomial) can count, for any $k$,
the number of $k$-component circuit partitions of the graph.
In this paper we define the interlace polynomial;
treat its connection to the Martin polynomial and
the Kauffman brackets of a link diagram;
and look into other apparently unrelated properties
of the polynomial, some of them only conjectural.
Basic questions about the polynomial remain open, but an
understanding is developing.
The connection with the Kauffman brackets was clarified, and that with
the Martin polynomial discovered, since the time of our first
publication~\cite{ABS00}.

\subsection{Looking backward}
The interlace graph polynomial has antecedents,
notably in the work of Bouchet. 
Specifically, starting from the study of
Euler circuits of undirected
4-regular graphs and the transformations of Kotzig~\cite{Kot68},
\cite{Bouchet87-isotropic-systems,
Bouchet88-graphic-presentations,
Bouchet91,
Bouchet-Tutte,
Bouchet-mm3}
include the pivot for graphs given by our \eqref{def:pivot},
and define a generalized Tutte-Martin polynomial on isotropic systems
(and on multimatroids),
of which as Bouchet~\cite{Bouchet-personal-comm00} pointed out
our interlace polynomial
may be viewed as a special case.
This is made explicit by Aigner and van der Holst~\cite{AignerHolst02}:
If a graph $G$ has vertex set $V$ and adjacency matrix~$A$,
$I$ is the identity matrix of the same dimension as~$A$,
and $\cal{S}_G = (V, \cal{L}_G)$ is the isotropic system
given by the row space $\cal{L}_G$ of $(A | I)$,
then our polynomial 
is equal to the Martin polynomial $m(\cal{S}_G;x)$.

\subsection{Looking forward}
In contrast to the reduction formula of Theorem~\ref{qdef},
the expression of the interlace polynomial as the Martin polynomial
of an isotropic system gives an
explicit expansion, as in Corollary~1 of \cite{AignerHolst02}.
Such an expansion was also produced,
without reference to isotropic systems, in \cite{2poly},
where a two-variable generalization of the interlace polynomial
is shown to have a formal similarity to the Tutte polynomial.

\subsection{Back to the present}
In this paper we introduce the one-variable interlace polynomial in the
circuit-enumeration context that led us to discover it,
and consider
extremal properties of the polynomial,
its values for various special graphs,
and evaluations which relate
to basic graph properties such as the component and independence numbers.

\subsection{Outline}
The paper is structured as follows.
In Section~\ref{sec:preliminaries} we define
an interlace graph $\HH$ for an Euler circuit $\ckt$
of a \tito directed graph $D$.
Section~\ref{sec:transposition}
defines a pivoting function on any undirected graph;
applied to an interlace graph $\HH$ the pivot is
consistent with transpositions on $C$, in a sense
specified by Lemma~\ref{lem:transpose}.
In Section~\ref{sec:circuits} this consistency is exploited to show
how $\HH$ determines the number of Euler circuits of~$D$,
giving a sort of preview of Sections \ref{sec:property} and~\ref{sec:ipoly}.

Section~\ref{sec:property} contains essential identities for the graph
pivot operator, and these are used in
Section~\ref{sec:ipoly}
to prove that the interlace graph polynomial $q$
is well-defined on any undirected graph~$G$,
not just on interlace graphs.

In the interlace-graph setting,
Section~\ref{sec:partitions} connects
$q(H)$ with
the circuit partition polynomial and
the Martin polynomial $m(D)$.
There is a further connection to the Kauffman brackets of a link diagram.

Section~\ref{sec:natural} computes $q(G)$ for some natural graphs,
and some general observations about the interlace polynomial are
drawn in Section~\ref{sec:facts}.
``Substituted'' and ``rotated'' graphs are considered in
Section~\ref{lesssimple}.
Calculating $q(G)$ on such graphs is
of mild interest in itself, and is essential machinery for
the computation in
Section~\ref{extremal} of extremal properties of the
graph polynomial.
Section~\ref{sec:open} considers some open problems.

\section{Interlacings} 
\label{sec:preliminaries}   

Our approach to counting circuits is based on ``interlacings'',
as in Read and Rosenstiehl~\cite{ReRo76}.
\begin{definition}
Given an Euler circuit $\ckt$ of a \tito directed graph $D$,
two vertices $\uu$ and $\vv$ of $D$ are \emph{interlaced}
if $\ckt$ visits them in the sequence
\ldots $\uu$ \ldots $\vv$ \ldots $\uu$ \ldots $\vv$ \ldots.
The \emph{interlace graph} $\HH=\ilace(\ckt)$ corresponding to $\ckt$
has the same vertex set as~$D$, with an edge $\uv$ in $\HH$ if $\uu$ and $\vv$
are interlaced in~$\ckt$.
\end{definition}
The class of graphs arising as interlace graphs coincides
with the
class of \emph{circle graphs},
as in~\cite{Golumbic,deFraysseix,Gabor,Spinrad}.
Trees constitute a particularly tractable set of interlace graphs.

\begin{figure}[bhtp]
\centerline{
\psfrag{u}[cc][cc]{$\mathbf{\uu}$}
\psfrag{v}[cc][cc]{$\mathbf{\vv}$}
\psfrag{1}[cc][cc]{$1$}
\psfrag{2}[Bc][Bc]{$2$}
\psfrag{3}[cc][cc]{$3$}
\psfrag{4}[cc][cc]{$4$}
\includegraphics[width=2.5in]{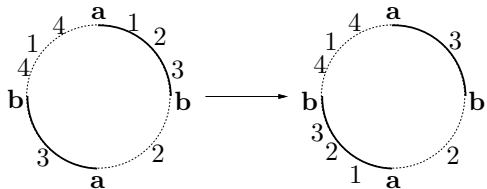}
}
\caption{In a circuit $\ckt$, transposing on $\uv$
toggles the interlacing relationships of some other pairs of vertices;
here 1 and 2, 1 and 3, and 2 and 3.
Vertices interlaced with $\uu$ are interlaced with $\vv$
after transposition, and vice-versa.
Also, the ``resolution'' at $a$ is reversed by transposition:
$\ckt$ follows edge $4a$ with $a1$ and $2a$ with~$a3$,
while $\ckt^{ab}$ follows $4a$ with $a3$ and $2a$ with~$a1$.
}
\label{fig:transpose2}   
\end{figure}

\begin{definition}
Given an Euler circuit $\ckt$ with $\uu$ and $\vv$ interlaced,
a \emph{transposition} on the pair $\uv$ is
the circuit $\ckt^{\uv}$ resulting from exchanging
one of the edge sequences from $\uu$ to $\vv$
with the other.
\end{definition}
See Figure~\ref{fig:transpose2} for an example.
Note that if $C$ traverses edge $e_1$ into $\uu$ and $e_2$ out on
one visit, and $e_3$ in and $e_4$ out on the other visit,
then $\ckt^{\uv}$ reverses this ``resolution'',
following $e_1$ with $e_4$, and $e_2$ with $e_3$.

The result stated below as Lemma~\ref{ukkpev} is due to
Ukkonen and Pevzner~\cite{Ukkonen,Pevzner}.  Their result is set in
a more general context, allowing ``rotations'' and ``3-way repeats'';
the relation between their context and \tito digraphs is described in detail
in Section~1.4 of~\cite{ABCS00}.  In the context of \tito digraphs,
we provide a simple proof of Lemma~\ref{ukkpev}, starting with the
case of a circuit having no interlaced pair.

\begin{lemma}
\label{ukkpev0}   
If a \tito digraph $D$ has an Euler circuit $C$ with no interlaced
pair, then $D$ has only one Euler circuit.
\end{lemma}

\begin{proof}
Use induction on the order $n$ of $D$.  For $n=1$ the only case is
$D$ consisting of two directed loops on a common vertex, for which there
is a single Euler circuit.  For $n>1$ the given circuit $C$ visits
the vertices in a sequence $s_1s_2\cdots s_{2n}$ in which
some vertex $a$ must appear twice consecutively.
(The two occurrences of each vertex can be viewed as a pair of
matched parentheses, and the absence of interlacings is equivalent
to the parentheses being balanced.  A vertex $a$ we seek corresponds
to an innermost pair of parentheses.)
At $a$ there is a loop; contracting out the loop, and then the
vertex $a$, gives a graph $D'$ with the same number of Euler circuits
as $D$. $D'$ also has no interlaced pair, and is of order $n-1$,
so by induction $D'$ has only one Euler circuit, and thus so does~$D$.
\end{proof}

\begin{lemma}
\label{ukkpev}   
For a \tito digraph $D$, all Euler circuits form a single orbit
under transpositions on interlaced pairs.
\end{lemma}

\begin{proof}
Induction on the order $n$ of~$D$.
The unique graph $D$ for the base case $n=1$ has only one Euler
circuit, so the assertion is trivial.
Otherwise, let two circuits $C$, $C'$ of $D$ be given.
If there is a vertex $a$ which the two circuits resolve alike
--- each of $a$'s two in-edges is followed by the same out-edge in
$C'$ as it is in $C$ ---
then use this resolution to contract $a$ out of $D$, $C$, and $C'$.
In the contracted case, by induction, one circuit can be converted
to the other by a sequence of transpositions, and this sequence
lifts to the original case.

If there is no vertex $a$ resolved alike in $C$ and $C'$, then $C$
and $C'$ are different, so by Lemma~\ref{ukkpev0}
$D$ must have had some interlaced pair~$ab$.
Transposing $C'$ on $ab$ gives a circuit $C''$ which resolves vertex
$a$ oppositely to~$C'$, and therefore the same as~$C$.
By the previous case, $C$ can be converted to~$C''$,
and in turn to~$C'$, by a sequence of transpositions.
\end{proof}

Lemma~\ref{ukkpev0} says that
if $\ilace(\ckt)$ has no edges,
then the Euler circuit is unique.
Also, if $\ckt$ has only one interlaced pair~$ab$,
transposing on it gives a circuit $\ckt'$ in which $ab$ is
still the sole interlaced pair,
so these two circuits constitute the full orbit;
that is, if $\ilace(\ckt)$ has one edge,
then there are exactly two Euler circuits.
Such observations,
studied in~\cite{ABCS00},
suggest a question: does the interlace graph
determine the number of Euler circuits?
We shall show that it does.
(In \cite{ABCS00} a more elementary fact is exploited,
that if a pairing generates $k$ Euler circuits,
it has at most $k-1$ interlaced pairs,
which define a sort of skeleton.
Thus the enumeration of pairings with $k$ Euler
circuits is reduced to an enumeration over a finite number of skeleta,
and a use of generating functions to count pairings with
a given skeleton.
While closely related, the interlace polynomial
does not simplify the enumeration.)

\section{The pivot operator on graphs}
\label{sec:transposition}   

We now define a pivot operator $G \mapsto G^{\uv}$ on any graph $G$.
It is connected with interlacings by Lemma~\ref{lem:transpose}:
if $\HH= \ilace(\ckt)$, and vertices
$\uu$ and $\vv$ are interlaced in $\ckt$,
then the circuit transposition $\ckt \mapsto \ckt^{\uv}$
and the pivot operator $\HH \mapsto \HH^{\uv}$ commute,
modulo a vertex relabelling.

\begin{figure}[bhtp]
\centerline{
\psfrag{u}[Br][Br]{$\mathbf{\uu}$}
\psfrag{v}[Br][Br]{$\mathbf{\vv}$}
\psfrag{toggle}[Bl][Bl]{\textsf{toggle}}
\includegraphics[width=3in]{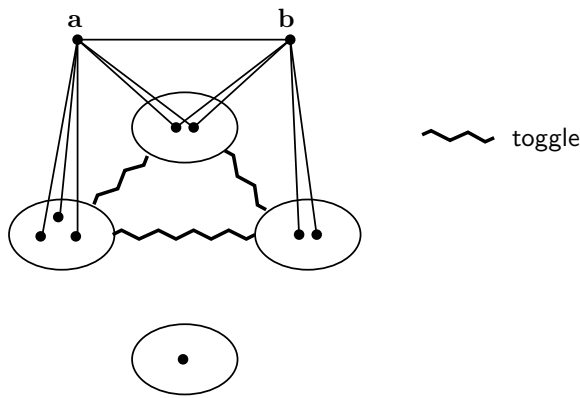}
}
\caption{Pivoting a graph $G$ on a pair $\uv$
(always an edge, in practice).
Vertices are divided into 4 classes:
those joined to $\uu$; to $\vv$; to both; or to neither.
If vertices $x$ and $y$ lie in different classes among the first
three classes named,
then if $xy$ was an edge in $G$, it becomes a non-edge in $G^{\uv}$,
and similarly if $xy$ was a non-edge, it becomes an edge.
}
\label{fig:transpose1}    
\end{figure}

Figure~\ref{fig:transpose1} illustrates the following definition.
\begin{definition}[Pivot]
\label{def:transpose}   
\label{def:pivot}   
Let $G$ be any undirected graph,
and $\uv$ a pair of distinct vertices in $G$.
Construct the \emph{pivot} $G^{\uv}$ as follows.
Partition the vertices other than $\uu$ and $\vv$ into four classes:
\begin{enumerate}
\item vertices adjacent to $\uu$ alone;
\item vertices adjacent to $\vv$ alone;
\item vertices adjacent to both $\uu$ and $\vv$;
and
\item vertices adjacent to neither $\uu$ nor $\vv$.
\end{enumerate}
Begin by setting $G^{\uv} = G$.
For any vertex pair $xy$ where $x$ is in one of the classes~(1--3)
and $y$ is in a different class~(1--3),
``toggle'' the pair $xy$:
if it was an edge, make it a non-edge,
and if it was a non-edge, make it an edge.
\end{definition}
Note that $G^{\uv}=G^{\vu}$.
Although the pivot operation is defined for any vertex pair $\uv$,
we shall only ever pivot about an \emph{edge}~$\uv$.

We will need one other, trivial, operation.
\begin{definition}[label-swapping]
\label{def:swap}   
For any directed or undirected labelled graph $G$
(or circuit $\ckt$),
and any pair of distinct vertices $\uu$ and $\vv$,
let $G_{\uv}$ (or $\ckt_{\uv}$)
denote the same object with the labels of $\uu$ and $\vv$ swapped.
\end{definition}

The pivot function defined above satisfies the following lemma.
\begin{lemma}
\label{lem:transpose}   
If $\HH = \ilace(\ckt)$,
and $\HH$ has an edge $\uv$ (\ie, $\uu$ and $\vv$ are interlaced in~$\ckt$),
then $\HH^{\uv} = (\ilace(\ckt^{\uv}))_{\uv}$.
\end{lemma}

\begin{proof}
A transposition on an Euler circuit $\ckt$ may be understood through
Figure~\ref{fig:transpose2}.
Vertices are written around the circumference of a circle in the
order they are visited by~$\ckt$,
and vertices $\uu$ and $\vv$ are interlaced if
the chord joining the two occurrences of $\uu$ intersects
the chord joining the two occurrences of~$\vv$.

First, note that after transposition, $\uu$ and $\vv$ remain interlaced.
Next, look at single vertices other than $\uu$ and $\vv$.
A vertex like $4$, interlaced with neither $\uu$ nor $\vv$ before transposition,
remains interlaced with neither.
A vertex like $3$, interlaced with both $\uu$ and $\vv$,
remains interlaced with both.
A vertex like $1$, interlaced with just $\uu$ in $\ckt$,
becomes interlaced with just $\vv$ in $\ckt^{\uv}$,
and therefore in $(\ckt^{\uv})_{\uv}$ is
interlaced with just $\uu$, as it was in $\ckt$.
After label-swapping, then,
incidences between $\uu$ and $\vv$ and other vertices are
unaffected by transposition.

Finally, consider incidences between pairs of vertices other than $\uu$ and~$\vv$.
In the example of Figure~\ref{fig:transpose2},
$4$ remains uninterlaced with $1$ and $3$,
but $1$ and $3$ become interlaced,
matching Definition~\ref{def:transpose}.
All other cases
can be worked out with similar figures,
and all accord with Definition~\ref{def:transpose};
this completes the lemma's proof.
\end{proof}

\section{Counting circuits in \tito digraphs}
\label{sec:circuits}   

The results in this section are a special case of those presented
in sections \ref{sec:property}--\ref{sec:partitions},
in particular Theorems \ref{qdef} and \ref{thm:q&r}.
However, the special case is relatively concise, and reveals
the interlace polynomial's relationship with Euler circuits
(its motivation) and the nature of the reduction it obeys.

\begin{definition}
\label{def:dk1}   
For any directed graph $D$, let $\numeuler(D)$ be the number of
Euler circuits of $D$, and, more generally, let $r_k(D)$ be the
number of partitions of the edges of $D$ into $k$ circuits.
\end{definition}

The fact that all Euler circuits of a \tito digraph
can be generated through transpositions
starting with any Euler circuit $\ckt$,
and transpositions are mirrored by
pivots of the interlace graph,
will mean that the number of Euler circuits of $D$
can be computed from the
interlace graph of any Euler circuit $\ckt$ of~$D$.
This is the basis of the following theorem.

\begin{theorem}
\label{thm:circuits} 
There exists a function $q_1$, from the set of interlace graphs
to the integers,
such that for any \tito digraph $D$ with Euler circuit $\ckt$,
the number of Euler circuits of $D$ is equal to $q_1(\ilace(\ckt))$.
Moreover, $q_1$ is uniquely defined by the following recursion:
\[
q_1(\HH) =
\begin{cases}
 q_1(\HH-\uu) + q_1(\HH^{\uv}-\vv) & \text{if $\uv$ is an edge of $\HH$} \\
 1 & \text{if $\HH$ has no edges.}
\end{cases}
\]
\end{theorem}

\begin{proof}
In the graph $D$, vertex $\uu$ has two in-edges, $e_1$ and $e_2$,
and two out-edges, $e_3$ and $e_4$.
(Since we allow loops, the in-edges are not
necessarily distinct from the out-edges.)
The set of all Euler circuits in $D$ can be partitioned into
those that follow $e_1$ with $e_3$
(and therefore follow $e_2$ with $e_4$),
and those that follow $e_1$ with $e_4$
(and therefore follow $e_2$ with $e_3$).
(We do not presume that either class is non-empty.)

It follows that
$\numeuler(D) = \numeuler(D') + \numeuler(D'')$.
Here $D'$ is derived from $D$ by merging edge $e_1$ with $e_3$
(replacing the two with a single new edge going from the tail of $e_1$
to the head of $e_3$),
merging $e_2$ with $e_4$,
and deleting the now-isolated vertex~$\uu$;
and $D''$ is similarly formed by merging $e_1$ with $e_4$,
merging $e_2$ with $e_3$, and deleting~$\uu$.

Define an analogous deletion operator on Euler circuits,
so that $\ckt-\uu$ means merging the edge preceding and the edge following
each occurrence of $\uu$, and removing~$\uu$.
Note that an Euler circuit $\ckt$ contains more information than its
parent graph $D$, and in particular $D=D(\ckt)$ is determined by~$\ckt$.
In the present case, we have $D=D(\ckt)$;
letting $\ckt' = \ckt-\uu$ and $\ckt'' = \ckt^{\uv}-\uu$,
we also have $D' = D(\ckt')$ and
(because the transposition $\ckt^{\uv}$ switches the allegiances
of the two edges into $\uu$ with the two out)
$D'' = D(\ckt'')$.
Defining $\numeuler(\ckt) = \numeuler(D(\ckt))$,
it is immediate that
$\numeuler(\ckt) = \numeuler(\ckt-\uu) + \numeuler(\ckt^{\uv}-\uu)$.

The interlace graphs of these circuits are given by
$\HH' = \ilace(\ckt-\uu) = \ilace(\ckt)-\uu$,
and (by Lemma~\ref{lem:transpose})
$\HH'' = \ilace(\ckt^{\uv}-\uu) = (\ilace(\ckt))^{\uv} - \vv$.
(In the last expression $\vv$ rather than $\uu$ is deleted because
of label-swapping; see Lemma~\ref{lem:transpose}.)
Since $\ckt-\uu$ and $\ckt^{\uv}-\uu$ each have one vertex less than $\ckt$,
suppose inductively that their Euler circuits are
counted as per the Theorem, that is,
$\numeuler(D') = q_1(\HH')$ and
$\numeuler(D'') = q_1(\HH'')$.
Then
\begin{align}
\numeuler(D) &= \numeuler(D') + \numeuler(D'') \notag
\\ &= q_1(\HH-\uu) + q_1(\HH^{\uv}-\vv)
 \label{eq:point} 
\\ &= q_1(\HH), \notag
\end{align}
the last equality following from the recursive definition of $q_1(\HH)$.
This shows both the independence of the quantity $q_1(\HH)$ from the
first pivot $\uv$ used in its recursive calculation,
and the desired equality $\numeuler(D) = q_1(\HH)$.
\end{proof}

The key point of this lies in equation~\eqref{eq:point},
showing the connection between
two in--out ``resolutions'' at a vertex $a$
and the interlace graph pivoted and not pivoted on an edge $ab$.

\section{Pivoting about an edge}
\label{sec:property}   

Consider the pivoting operation of Definition~\ref{def:transpose},
and recall that $\Gab=G^{ba}$.
Also note that pivoting is an involution:
for any pair $ab$, $G^{(ab)(ab)}=G$,
where the notation $G^{(ab)(cd)}$
represents an iterated pivot $(\Gab)^{cd}$.
We shall need the following properties of pivoting,
which are considerably more subtle.

\begin{lemma}
\label{lem1}   
Let $a$, $b$, and $c$ be distinct vertices of a graph~$G$ with $ab,ac\in E(G)$.
Then
\begin{itemize}
\item[(i)] $G^{(ab)(ac)(ab)}=G_{bc}$, and
\item[(ii)] $G^{(ab)(ac)}=(G^{ac})_{bc}$.
\end{itemize}
\end{lemma}

Two remarks before the proof.
The identities need not hold unless $ab$ and $ac$ are both edges.
And, the identities are needed to prove Lemma~\ref{lem:pivot},
which in turn is needed for Theorem~\ref{qdef}.

\begin{proof}
Focus first on identity~(i), \emph{in the case where} $bc \notin E(G)$.
Proof of (i) consists of checking three types of edges.

First, inspection of Definition~\ref{def:pivot}
shows that the presence or absence in $G^{(ab)(ac)(bc)}$
of any edge in $\{a,b,c\} \times \{a,b,c\}$ is determined
by $G \restricted{\{a,b,c\}}$,
the subgraph of $G$ induced by $a$, $b$, and~$c$.
Checking (i) here merely requires carrying out the pivot operations
on the three-vertex graph with edges $ab$ and $ac$.
(This and the subsequent checks described can certainly
be carried out by hand, but we did it by computer instead,
with a small amount of manual spot-checking.)

Second, for any vertex $u \notin \{a,b,c\}$,
the presence or absence in $G^{(ab)(ac)(bc)}$ of
any edge in $\{u\} \times \{a,b,c\}$ is determined
by $G \restricted{\{a,b,c,u\}}$.
Checking (i) here requires checking 8 cases, one for each
setting (present/absent) of the three edges $ua$, $ub$, and $uc$.
They may be checked ``in parallel'' by computing
$G^{(ab)(ac)(bc)}$ for a single graph $G$ with vertices
$a$, $b$, $c$, and 8 more vertices belonging to different
equivalence classes as defined by adjacency to $a$, $b$, and~$c$.
Let us define $G$ not to include any edges except those incident
on $\{a,b,c\}$, as such edges are irrelevant to this calculation.

Finally, for any two vertices $u,v \notin \{a,b,c\}$,
the number of times the presence/absence of the edge $uv$ is
toggled depends on the number of times $u$ and $v$ lie in different
classes (1--3) of Definition~\ref{def:pivot}.
But the class of $u$ (respectively $v$) is determined
by $G \restricted{\{a,b,c,u\}}$.
Note then that the number of times $uv$ is toggled is independent
of the initial presence or absence of $uv$ or any other edge
not incident on $\{a,b,c\}$.
Thus here too it suffices to compute
$G^{(ab)(ac)(bc)}$ for the same 11-vertex graph $G$
as above.

In short, in the case $bc \notin E(G)$,
identity (i) can be verified by computing
$G^{(ab)(ac)(bc)}$ for a single 11-vertex graph
with vertices $a$, $b$, and $c$, and 8 other vertices representing
the 8 equivalence classes defined by adjacency to $a$, $b$, and~$c$.

\emph{In the case} $bc \in E(G)$ the identity is checked the same way,
using the same 11-vertex graph as above but with the edge $bc$ added.

Identity (ii) could be proved using the same pair of graphs,
but in fact, (i) and (ii) are equivalent.
Since $G \mapsto \Gab$ is an involution,
$G^{(ab)(ac)(ab)}=G_{bc}$ iff
$G^{(ab)(ac)}=(G_{bc})^{ab}$.
But $(G_{bc})^{ab}=(G^{ac})_{bc}$, so we are done.
\end{proof}

\begin{remark}
\label{pivot-connected}   
If $G$ is connected, then for any edge $ab$ of $G$, $\Gab$ is connected.
\end{remark}
\begin{proof}
Edge relations change only between vertices ``distinguished'' by $a$ and $b$,
but in $\Gab$ these vertices remain joined to $a$ or $b$ (or both),
and $a$ and $b$ remain joined, so $\Gab$ is connected.
\end{proof}

\section{The interlace polynomial}
\label{sec:ipoly}   

We will now show that the integer $q_1(G)$
defined by Theorem~\ref{thm:circuits}
can be generalized to a one-variable graph \emph{polynomial}
$q(G;x)$ which is defined for \emph{all} graphs $G$,
not just interlace graphs, and which for interlace graphs $G$
satisfies $q(G;1) = q_1(G)$.

We shall call $q(G)$ the \emph{interlace polynomial} of $G$,
and shall write $q(G;x)$, $q(G)(x)$, and $q_G(x)$ interchangeably.
Later we will show that when $\HH$ is an interlace graph of $D$,
$q_\HH(x)$ can be used to count partitions of $D$ into circuits.

Historically, we began with the graph function $q_1$,
whose well-definedness --- on interlace graphs ---
followed from its construction.
Generalizing to a graph \emph{polynomial} $q$ was a leap of faith,
but computer experiments quickly convinced us that
the polynomial of Theorem~\ref{qdef} was indeed
well-defined, on \emph{arbitrary} graphs ---
that the sequence of pivot elements chosen was irrelevant.
We expected that to prove it we would need to prove something
like Lemma~\ref{lem:pivot}, \ie, that we would
have to prove the equivalence of any two \emph{first}
pivot elements, and that this could be boiled down to the
equivalence of $q^{ab}$ to $q^{ba}$,
and that of $q^{ab}$ to $q^{ac}$ (or perhaps to $q^{ca}$).
Then, we guessed that showing for example that
$q^{ab}(G) = q^{ac}(G)$ would mean showing that
following the $ab$ pivot with one on $ac$ or $ca$
gives the same set of \emph{graphs} as
following the $ac$ (or $ca$) pivot with one on $ab$ or $ba$.
Computer experiments then indicated the correct identities,
leading to the rather obscure but easy to prove identities
of Lemmas \ref{lem1} and~\ref{lem:pivot}.

\subsection{Definition and well-definedness}

As usual, let us write ${\cal G}$ for the class of finite undirected graphs
having no loops nor multiple edges.

\begin{theorem}[Interlace polynomial]
\label{qdef}   
There is a unique map $q:{\cal G} \to \mathbb{Z}[x]$, $G \mapsto q(G)$, such
that the following two conditions hold.
\begin{itemize}
\item[(i)] If $G$ contains an edge $ab$ then
\begin{align}
\label{eq:qdef}   
q(G)=q(G-a)+q\left(\Gab-b\right) .
\end{align}
\item[(ii)] On $E_n$, the edgeless graph with $n$ vertices, $q(E_n)=x^n$.
\end{itemize}
\end{theorem}
The only issue in proving the theorem is to show that $q$
is well-defined, independent of the pivot order.

Our introductory assertion is an immediate corollary of the theorem.
\begin{corollary}
If $G$ is an interlace graph,
$q_1(G) = q(G;1)$.
\end{corollary}

The following Lemma is used to prove Theorem~\ref{qdef}.
\begin{lemma}
\label{lem:pivot}   
Suppose that, for all vertices $a$ and $b$ of $G$,
$q(G-a)$ and $q(\Gab-b)$ are both well defined, and
let $q^{ab}(G) \eqdef q(G-a)+q(\Gab-b)$.
\begin{itemize}
\item[(i)]
If $ab$ is an edge of~$G$ then
$$
q^{ab}(G)=q^{ba}(G).
$$
\item[(ii)]
If $ab$ and $ac$ are edges of~$G$ then
$$
q^{ab}(G)=q^{ac}(G).
$$
\end{itemize}
\end{lemma}

\begin{proof}
Both parts of the proof will be by induction on the order of~$G$.

(i) If $d(a)=1$ or $d(b)=1$
then $\Gab=G$, so
$$
q^{ab}(G) \; = \; q(G-a)+q(\Gab-b) \; = \; q^{ba}(G).
$$
Hence, we may assume that $G$ contains edges $ac$ and $bd$,
distinct from $ab$.
(We do not assume that $c$ and $d$ are distinct.)
Start with
\begin{align}
q^{ab}(G) &= q(G-a)+q(\Gab-b) . \notag \\
\intertext{By the inductive hypothesis, both the latter terms are
well defined, and we may apply any sequence of pivot edges we like.
Pivot the first term on $bd$,
and the second on $ac$, which is an edge of $\Gab$ as it was of $G$:}
&=
\{q(G-a-b)+q((G-a)^{bd}-d)\} +
 \notag \\ & \ind \{q(\Gab-b-a)+q((\Gab-b)^{ac}-c)\} \notag \\
&= q(G-a-b)+q(G^{bd}-a-d) +
 \label{ab-expand} 
 \\ & \ind
 q(\Gab-a-b) + q(G^{(ab)(ac)}-b-c) . \notag
\\ \intertext{%
Symmetrically,
swapping the roles of $a$ and $b$, and those of $c$ and $d$,
}
q^{ba}(G)
&= q(G-a-b)+q(G^{ac}-b-c)
 \label{ba-expand} 
 \\ & \ind
 + q(\Gab-a-b) + q(G^{(ab)(bd)}-a-d) . \notag
\end{align}
In these expansions of $q^{ab}(G)$ and $q^{ba}(G)$, the first
terms are identical, as are the third terms.
Also, by Lemma~\ref{lem1}~(ii),
the second term of \eqref{ab-expand} is equal to the fourth of
\eqref{ba-expand}, because
$G^{bd}-a-d \isom G^{(ab)(bd)}-a-d$, and
similarly the fourth term of \eqref{ab-expand} is equal to the
second of \eqref{ba-expand},
so $q^{ab}(G)=q^{ba}(G)$.

(ii) This identity is also an immediate consequence of Lemma~\ref{lem1},
and is proved similarly,
but here the pivot sequences needed do not respect symmetry.
\begin{align*}
q^{ab}(G) &= q(G-a)+q(\Gab-b) ; \\
\intertext{both these terms are well-defined by induction,
and pivoting the second on $ac$, }
&=
q(G-a) +\{q(\Gab-b-a) +
q((\Gab-b)^{ac}-c)\}\\
&= q(G-a)+q(\Gab-a-b)+q(G^{(ab)(ac)}-b-c) .
\end{align*}
And,
\begin{align*}
q^{ac}(G) &= q(G-a)+q(G^{ac}-c) ; \\
\intertext{pivoting the second term on $ba$
(\emph{not} on $ab$ as symmetry would suggest), }
&=
q(G-a) + \{ q(G^{ac}-c-b) +
q((G^{ac}-c)^{ba}-a) \} \\
&= q(G-a)+q(G^{ac}-b-c)+q(G^{(ac)(ab)}-a-c) .
\end{align*}
By Lemma~\ref{lem1}, $(\Gab-a-b)_{bc}=G^{(ac)(ab)}-a-c$
and $G^{(ab)(ac)}-b-c=G^{ac}-b-c$,
so $q^{ab}(G) = q^{ac}(G)$.
\end{proof}

Now we are ready to prove the main theorem. 

\begin{proof}[Proof of Theorem~\ref{qdef}]
If $G$ is a graph of order~$n$ \emph{containing an edge} $ab$
then the pivot-reduction formula \eqref{eq:qdef}
reduces the computation of $q(G)$
to the computation of~$q$ on two graphs of order $n-1$, namely
$G-a$ and $\Gab-b$.  Repeating this process, we reduce $q(G)$ to a
linear combination $q(G)=\sum_{k=1}^{n-1}a_kq(E_k)$, where each $a_k$
is a nonnegative integer.  By making use of the
boundary conditions $q(E_k)=x^k$, we find $q(G)$.
Hence, \emph{if there is a function} $q$ then it is unique.

However, it is not clear that \emph{there is} a function~$q$
satisfying (i) and~(ii).  To prove this, we have to show that
the same result arises irrespective of the pivot edges
chosen in the recursive reduction of $G$ to edgeless graphs.
The proof will be by induction on the order $n$ of $G$.
At $n=1$ there is nothing to prove, so assume that $n>1$.
By the inductive hypothesis, $q(G-a)$ and $q(\Gab - b)$
are both well defined, so we need only show that
$q(G-a) + q(\Gab-b)$ is independent of the
\emph{first} pivot edge $ab$ chosen.
(Note that $q^{ab} \eqdef q(G-a) + q(\Gab-b)$,
unlike $\Gab=G^{ba}$,
depends on the \textit{ordered} pair $(a,b)$, that is, on the
\textit{oriented} edge $ab$.)
To restate, we must show that for any two oriented edges $ab$ and $a'b'$,
$q(G-a) + q(\Gab-b) = q(G-a') + q(G^{a'b'}-b')$.

We consider two cases.
The first comes if $ab$ and $a'b'$ lie in two different
components of~$G$:
say $a,b \in G_1$ and $a',b' \in G_2$, with $G = G_1 \cup G_2$.
Pivoting $G$ on $ab$ does not affect the component $G_2$ at all,
since all $G_2$'s vertices fall into Definition~\ref{def:pivot}'s case~(ii).
Thus
\begin{align*}
q(G-a) + q(\Gab-b)
 &= q((G_1 -a) \cup G_2) + q((G_1^{ab}-b) \cup G_2) .
\\ \intertext{By the inductive hypothesis,
each of the latter terms is well defined,
and we may continue as we like.  Pivoting on $a'b'$,}
 &= q((G_1-a) \cup (G_2-a')) + q((G_1^{ab}-b) \cup (G_2^{a'b'}-b')) .
\end{align*}
By symmetry, this is clearly equal to pivoting first on $a'b'$
and then on $ab$; thus
$q(G-a) + q(\Gab-b) = q(G-a') + q(G^{a'b'}-b')$
as was required.
Note that isomorphism of $G$ and $G'$ is our only tool for proving
equality of $q(G)$ and $q(G')$, so our method was (and always will be)
to apply pivot-reduction to each of $q(G)$ and $q(G')$ until the graphs
comprising their sums are in one-to-one correspondence.

The more interesting case comes if $ab$ and $a'b'$
lie in the same component of $G$.
If so, $G$ contains a walk $a,b, \ldots, a', b'$.
Lemma~\ref{lem:pivot} shows that for any sub-walk $v_i, v_{i+1}, v_{i+2}$,
we have
$q^{v_i, v_{i+1}}(G) = q^{v_{i+1}, v_{i+2}}(G)$,
whether $v_{i+2}$ is distinct from $v_i$ (the usual case),
or $v_{i+2} = v_i$
(the case we use to re-orient a terminal edge if necessary).
Given this, $q^{ab}(G) = \cdots = q^{a'b'}(G)$ as desired.
\end{proof}

While we will not exploit the following mild generalization
of Theorem~\ref{qdef},
it demonstrates that the structure of polynomial multiplication
is not required:
the powers $x$, $x^2$, \ldots
could be replaced by indeterminates $x_1$, $x_2$, \ldots.

\begin{theorem}
There is a unique map $\qq:{\cal G} \to \mathbb{Z}[x_1, x_2, \ldots]$,
the $\mathbb{Z}$-module with basis $x_1, x_2, \ldots$,
such that the following two conditions hold.
\begin{itemize}
\item[(i)] If $G$ contains an edge $ab$ then
$$
\qq(G)=\qq(G-a)+\qq\left(\Gab-b\right).
$$
\item[(ii)] On the edgeless graph $E_n$, $\qq(E_n)=x_n$.
\end{itemize}
\end{theorem}
\begin{proof}
Identical to the proof of Theorem~\ref{qdef}.
\end{proof}

On the other hand, polynomial multiplication does have a natural role
for the interlace polynomial,
as shown by Remark~\ref{G1G2},
the first of several properties we now consider.

\subsection{Simple properties}
\label{simpleProperties}   
We should first remark that $q(G;x)$ does not appear
to be a specialization of the Tutte polynomial $T(G;x,y)$.
At least, for any tree $G$ of order $n$ the Tutte polynomial
is always $T(G) = x^{n-1}$, while the interlace polynomial varies
widely over trees.

With $G_1 = (V_1, E_1)$ and $G_2 = (V_2, E_2)$,
and $V_1 \cap V_2 = \emptyset$,
let $G = G_1 \cup G_2$ denote the disjoint union of the two graphs,
so $G = (V_1 \cup V_2, E_1 \cup E_2)$.
The following remark is trivial but frequently employed.

\begin{remark}
\label{G1G2}   
For any graphs $G_1$ and $G_2$ on disjoint vertex sets,
$q(G_1 \cup G_2) = q(G_1) q(G_2)$.
\end{remark}

For the Tutte polynomial, the corresponding equality
$T(G_1 \cup G_2) = T(G_1) T(G_2)$ holds even if $G_1$ and $G_2$
share a single vertex,
but the path $P_2$ suffices to demonstrate that this generalization
does not hold for the interlace polynomial.

\begin{remark}
\label{order}   
For any graph $G$ of order $n$, $q_G(2) = 2^n$.
\end{remark}
\begin{proof}
If $G=E_n$ the conclusion is immediate.  Otherwise it follows
from a trivial induction on $n$:
$q(G;2) = q(G-a;2)+q(\Gab-b;2) = 2^{n-1}+2^{n-1}=2^n$.
\end{proof}

\begin{remark}
\label{qab}   
For any graph $G$ and any edge $ab$,
$q(\Gab) = q(G)$.
\end{remark}

\begin{proof}
Pivoting $G$ on $ab$, $q(G) = q(G-a) + q(\Gab-b)$.
Pivoting $\Gab$ on $ba$,
$q(\Gab) = q(\Gab-b) + q(G^{(ab)(ba)}-a) = q(\Gab-b) + q(G-a)$.
\end{proof}

\begin{remark}
\label{indep}   
Let $G$ be a graph and let $H$ be any induced subgraph of $G$.
Then $\deg(q(G)) \geq \deg(q(H))$.
In particular, $\deg(q(G)) \geq \alpha(G)$
where $\alpha(G)$ is the independence number, \ie,
the maximum size of an independent set.
\end{remark}

\begin{proof}
For the main result,
all one has to show is that $\deg(q(G)) \geq \deg(q(G-a))$.
This holds because
if $a$ is an isolated vertex, then $q(G) = x q(G-a)$;
otherwise, for any edge $ab$,
$q(G) = q(G-a) + q(\Gab-b)$,
and the non-negativity of all terms implies $q(G) \geq q(G-a)$.

The corollary follows by choosing $H$ to be a maximum independent
set in~$G$.
\end{proof}

We do not know if there exist graphs for which
$\deg(q(G)) \gg \alpha(G)$.

\begin{remark}
\label{component}   
The degree of the lowest-degree term of $q(G)$ is precisely $k(G)$,
the number of components of~$G$.
\end{remark}
\begin{proof}
From Remark~\ref{G1G2}, it suffices to show that for any connected
graph $G$, $q(G)$ has a linear term and no constant term.
That the constant term of $q$ is 0 whenever $\order{G} \geq 1$
is trivial by induction.
We also prove by induction that $q(G)$ has a linear term.
Any connected graph $G$ has some non-cut vertex $a$.
Pivoting on any edge $ab$, $q(G)=q(G-a)+q(\Gab-b)$,
$q(G-a)$ is connected and is subject to the inductive hypothesis,
$q(\Gab-b)$ adds no constant term,
and so $q(G)$'s lowest-order term is linear.
\end{proof}

\def\K{\mu}

\begin{remark}
\label{matching}   
Let $\K(G)$ denote size of a maximum matching
(maximum set of independent edges)
in a graph $G$, and let $\deg(q)$ denote the degree of the
polynomial $q$.  If $G$ is a forest with $n$ vertices,
then $\deg(q(G))=n-\K(G)$.
\end{remark}

The straightforward inductive proof is omitted.

Taken together, Remarks \ref{order}, \ref{indep},
\ref{component}, and \ref{matching}
mean that from the interlace polynomial $q(G)$ we can read off the
graph's order,
an upper bound on its independence number,
its component number,
and (if $G$ is a forest) its matching number.

\section{Circuit partitions and the Martin and interlace polynomials}
\label{sec:partitions}   

In this section we relate the interlace polynomial to the
enumeration of circuit partitions of a \tito digraph.
Various circuit-counting polynomials were defined by
Martin~\cite{Martin77},
and the topic has been explored
by Martin~\cite{M2},
Las Vergnas~\cite{LV81,LV83,LV88},
Jaeger~\cite{J1},
Ellis-Monaghan~\cite{E1,E99},
and Bollob\'as~\cite{Boll02}.

If $\HH= \ilace(\ckt)$ is the interlace graph of an Euler circuit $\ckt$
of a \tito digraph~$D$ of order~$n$,
then by Theorem~\ref{thm:circuits},
$q_\HH(1)$ is the number of Euler circuits of~$D$.
Since every pairing of the in- and out-edges at each
vertex of $D$ leads to a \emph{circuit partition} of $D$,
there are $2^n$ such circuit partitions;
also,
it can be seen immediately from Theorem~\ref{qdef} that $q_\HH(2) = 2^n$.
In fact, the entire polynomial $q_\HH(x)$
can be interpreted in terms of decompositions of $D$ into circuits.

Let us define a somewhat unusual class of objects.
We call a digraph \emph{weakly Eulerian}
if it is an edge-disjoint union of oriented circuits
and ``free loops'',
that is loops without vertices.
Equivalently, a weakly Eulerian digraph is
a digraph in which every vertex has the same in-degree as
out-degree, together with some free loops.
These graphs may have multiple edges (and multiple loops),
all of which are distinguishable.

Definition~\ref{def:dk1}, of polynomial $r_k(D)$ counting the
partitions of a digraph $D$ into $k$ circuits, extends naturally
to weakly Eulerian digraphs.
Note that $\sum_k r_k(D) > 0$ if and only if $D$ is weakly Eulerian.

\begin{definition}
\label{def:decomp}   
For any weakly Eulerian digraph $D$,
the directed circuit partition polynomial is
$$
r(D)(x)=\sum_{k} r_k(D) x^k.
$$
\end{definition}

For example if $D$ consists of 3 loops on a single vertex,
$r_1(D) = 2$, $r_2(D) = 3$, and $r_3(D) = 1$.
As observed by Las Vergnas~\cite{LV83},
if $D$ consists of $m$ loops on a single vertex then
$r_k(D) = s(m,k)$, an unsigned Stirling number of the first kind
(see for example~\cite{Wilf} and~\cite{PolyaNotes}),
and $r(D)(x) = x(x+1) \cdots (x+m-1)$.
Note that $r_0(D)=0$ unless $D$ has no edges, in which case $r(D)=1$.

It is immediate that
if $D_1,\ldots,D_k$ are the
components of $D$ then $r(D) = \prod_{i=1}^k r(D_i)$.
To compute $r(D)$, choose any vertex $a$ of $D$,
and let $D_\sigma$ be the various resolutions at~$a$
(matchings of in-edges to out-edges,
followed by contraction of the vertex~$a$);
then $r(D)  = \sum_\sigma r(D_\sigma)$.

The circuit partition polynomial is simply a transformation of
the Martin polynomial.
Introduced by Martin in 1977~\cite{Martin77},
this polynomial was studied extensively by Las Vergnas~\cite{LV83},
who is responsible for its generalization to all Eulerian graphs.

\begin{definition}
\label{Martin}   
For any weakly Eulerian directed graph $D$, the Martin polynomial is
$$ m(D;x) = \frac1{x-1} r(D; x-1) . $$
\end{definition}

Recall that $q(G)$ is defined for any undirected graph $G$,
and we have just defined the circuit partition polynomial
$r(D)$ and the Martin polynomial $m(D)$ for any weakly Eulerian digraph~$D$.
In the event that $D$ is an Eulerian (connected) \tito digraph,
and $G$ an associated
interlace graph, the following theorem describes
an interesting connection between them.

\begin{theorem}
\label{thm:q&r}   
Let $D$ be an Eulerian
\tito digraph, $\ckt$ any Euler circuit of $D$,
and $\HH = \ilace(\ckt)$ its interlace graph.
Then
\begin{align}
q(\HH;x) = \frac1{x-1} r(D; x-1) = m(D;x) .
\label{d2a}   
\end{align}
\end{theorem}

\begin{proof}
The second equality is just Definition~\ref{Martin}.
For the first, we prove the transformation $x q(\HH;1+x) = r(D;x)$.
Since $\ckt$ conveys more information than $D$,
let us define $r(\ckt) = r(D)$.
First, suppose some $\uu$ and $\vv$ are interlaced in $\ckt$
(so $\uv$ is an edge in $\HH$).
Recall that, as in the proof of Theorem~\ref{thm:circuits},
$\ckt$ and $\ckt^{\uv}$ ``resolve'' the vertex $\uu$ in the
two different ways.
Thus,
\begin{align*}
r(D;x) &= r(C;x)
 \\ &= r(C-\uu;x) + r(C^{\uv}-\uu;x) \\
 \intertext{which by induction on the order of $D$ (and of $C$ and $\HH$)}
 &= x q(\HH-\uu; 1+x) + x q(\HH^{\uv}-\vv; 1+x)
 \\ &= x q(\HH;1+x).
\end{align*}

If there is no interlaced pair $\uu,\vv$ in $\ckt$, \ie, if
$\HH$ is an edgeless graph, then $\ckt$ must have a loop on some vertex,
and this loop may either be traversed as a separate circuit,
giving $xr(C-\uu;x)$, or may combine with the rest, giving $r(C-\uu;x)$;
in all, this gives $r(C;x) = (1+x) r(C-\uu;x)$, and induction
completes the proof.
\end{proof}

There is a trivial generalization to disconnected \tito graphs.
Such a graph $D$ consists of vertex-disjoint \tito components,
$D = D_1 \cup \cdots \cup D_n$.
Each $D_i$ has an Euler circuit $C_i$,
and if we \emph{define} an interlace graph $H = \bigcup H_i$ to be the
union of the components' interlace graphs,
then $r(D) = \prod r(D_i)$ and $q(H) = \prod q(H_i)$,
and it is immediate from Theorem~\ref{thm:q&r}
that $r(D) = xq(H;1+x)$.

Theorem~\ref{thm:q&r} has the following consequence for the
coefficients of the circuit-counting and interlace polynomials.

\begin{corollary}
\label{cor:q&r&m}   
Let $D$ be an Eulerian \tito digraph, $\ckt$ any Euler circuit of $D$,
and $\HH = \ilace(\ckt)$ its interlace graph.
If $D$ has $r_k$ partitions into $k$ circuits,
and $q(G;x) = \sum_k a_k x^k$,
then
\begin{align*}
r_k=\sum_{\ell} a_{\ell} \binom{\ell}{k-1}
\text{ \qquad and \qquad }
a_k = \sum_{\ell} r_{\ell+1} (-1)^{\ell-k} \binom{\ell}{k}  .
\end{align*}
\end{corollary}

\begin{proof}
The first part is immediate from the preceding theorem and
the binomial formula.
The second follows from the
inclusion-exclusion principle (see for example~\cite{PolyaNotes}),
in particular, from the fact that the matrix $M$ with entries
$m_{ij} = \binom{i}{j}$, $i$ and $j$ ranging from $0$ to $n$,
has as its inverse a similar matrix with entries
$(m^{-1})_{ij} = (-1)^{i+j} \binom{i}{j}$.
\end{proof}

Corollary~\ref{c-2} below makes another connection between
the Martin polynomial,
and the interlace polynomial of an interlace graph.
To introduce it we need a couple of simple preliminaries and
a result of Martin's.
Let $a$ be a vertex of a \tito digraph $D$ of order $n$,
with arcs $e_1$, $e_2$
entering $a$ (say from $u_1$ and $u_2$)
and arcs $f_1$, $f_2$ leaving $a$ (say to $v_1$ and $v_2$).
Let $D_a'$ and $D_a''$
be the digraphs obtained from $D$ by resolving $a$ in the two possible
ways by uniting an incoming edge with an outgoing edge. For example,
$D_a'$ is obtained from $D-a$ by adding to it an arc from $u_1$ to
$v_1$ and an arc from $u_2$ to $v_2$, and $D_a''$ is obtained from
$D-a$ by adding an arc from $u_1$ to $v_2$ and one from $u_2$ to
$v_1$. If there is a loop at $a$, so that $u_1=v_1=a$, say, then
$D_a'$ consists of a \tito digraph with  $n-1$ vertices and a
vertexless loop, and in $D_a''$ all arcs incident with $a$ are replaced
by a single arc. Note that
\begin{equation}
\label{rres}   
r(D)=r(D_a')+r(D_a'').
\end{equation}

An {\em anti-circuit} is a circuit whose consecutive
arcs have opposite orientations.
Thus a cycle (cycle meaning circuit with no repeated vertex)
is an anti-circuit if every other vertex has
in-degree 0 and out-degree~2, the remaining vertices having
in-degree 2 and out-degree~0.
Every \tito digraph $D$ has a unique decomposition into anti-circuits;
we write $a(D)$ for the number of anti-circuits in this decomposition.
Martin made the following connection between
anti-circuits and the value of $m(D)$ at~$-1$;
our proof is much simpler than that in Martin's thesis~\cite{Martin77}.

\begin{theorem}[Martin]
\label{r-2}   
Let $D$ be a \tito digraph with $n$ vertices. Then
\begin{equation*}
-2 m(D; -1) = r(D; -2)=(-1)^{n+a(D)} 2^{a(D)}.
\end{equation*}
\end{theorem}

\begin{proof}
The first equality is immediate from~\eqref{d2a}.  For the second,
we apply induction on $n$. For $n=1$ our graph $D$ consists of two
loops on a vertex, so
$a(D)=1$, $r(D; x)=x+x^2$ and so $r(D; -2)=2=(-1)^{1+1} 2^1$,
as claimed. Suppose then that $n > 1$
and the theorem holds for smaller values of~$n$. Let $D_a'$ and
$D_a''$ be as in~\eqref{rres}. If $a$ is contained in two alternating
circuits then $a(D_a')=a(D_a'')=a(D)-1$, so
\begin{align*}
r(D; -2)&=r(D_a'; -2)+r(D_a''; -2)\\
&= 2 \cdot (-1)^{(n-1)+(a(D)-1)} 2^{a(D)-1}
 =(-1)^{n+a(D)} 2^{a(D)}.
\end{align*}
If on the other hand $a$ is contained in a single alternating
oriented circuit, then $\{a(D_a'), a(D_a'')\}=\{a(D), a(D)+1\}$, so
\begin{align*}
r(D; -2)
&= (-1)^{(n-1)+a(D)} 2^{a(D)} + (-1)^{(n-1)+(a(D)+1)} 2^{a(D)+1}
 = (-1)^{n+a(D)} 2^{a(D)},
\end{align*}
completing the proof.
\end{proof}

\begin{corollary}
\label{c-2}   
Let $H$ be an interlace graph. Then
\begin{equation*}
\label{circleminus}   
q_H(-1)=(-1)^{|H|+1}(-2)^k
\end{equation*}
for some non-negative integer $k$.
\end{corollary}

This led us to conjecture in \cite{ABS00} that for all undirected
graphs $G$, not just interlace graphs, $q(G;-1) = (-1)^{|G|} (-2)^k$
for some value~$k$.
This has recently been proved
by Balister, Bollob\'as, Cutler, and Pebody~\cite{BBCP00},
which for the interlace polynomial extends a result
of Rosenstiehl and Read on the bicycle dimension
of a graphical matroid~\cite{RR}.
The conjecture's implications may be at least as interesting
as the conjecture itself.
If $ab$ is an edge of $H$ and $G=H-a$, then $q(H;-1)=q(G;-1)+q(\Hab-b;-1)$.
For all three to be powers of~2, it must be that
$q(H;-1) \in \{\frac12 q(G;-1), -q(G;-1), 2 q(G;-1)\}$.
That is, adding a vertex to any graph $G$ in any way multiplies $q(-1)$
by a factor of $1/2$, $-1$, or $2$;
symmetrically, deleting any vertex from any graph $H$ multiplies $q(-1)$
by a factor of $1/2$, $-1$, or $2$.

In conclusion, let us point out a connection between the
interlace polynomial and the Kauffman bracket.
Given an alternating link diagram $L$ with set of crossings $V$,
let $D$ be the \tito digraph with vertex set $V$ whose edges
are the strands of $L$,
with each strand oriented from its over-crossing to its under-crossing.
Writing $[L](A,B,d)$ for the Kauffman square bracket,
it follows from results of Martin~\cite{martin78} and
Kauffman~\cite{kauffman87,kauffman88} that
\begin{align}
[L](1,1,x) &= \frac{1}{x}r_D(x) = m_D(x+1) = q_\HH(x+1) ;
\end{align}
these identities are also immediate consequences of the definitions.

An anonymous referee points out another connection, in the same setting.
If $G$ is a planar graph and $G_m$ is its medial graph,
oriented by going counterclockwise around the ``black'' faces,
then the Tutte polynomial is related to the Martin polynomial by
$T(G;x,x) = m(G_m;x)$
\cite{Martin77,LV81},
and in turn to the interlace polynomial through Theorem~\ref{thm:q&r}.

\section{Interlace polynomials of some simple graphs}
\label{sec:natural}   

The following results are all quite simple and so the proofs are omitted
for brevity.

\begin{example}
\label{easygraphs}   
The interlace polynomials of edgeless graphs, complete graphs,
stars, and complete bipartite graphs are given by
\begin{align*}
q(E_n) &= x^n           & \text{for} \ n & \geq 0
\\
q(K_n) &= 2^{n-1}x      & \text{for} \ n & \geq 1
\\
q(K_{1n})
           &= 2x + x^2 + x^3 + \cdots + x^n
                        & \text{for} \ n & \geq 2
\\
q(K_{mn}) &= (1+ \cdots + x^{m-1}) (1+ \cdots + x^{n-1}) + x^m + x^n - 1
                        & \text{for} \ m,n & \geq 1 .
\end{align*}
\end{example}

\begin{example}
For $n\ge 2$, the interlace polynomial of the path $P_n$
(with $n+1$ vertices and $n$ edges)
satisfies
\begin{align}
\label{eqpath1}   
q(P_n) &= q(P_{n-1})+xq(P_{n-2}),
\\
\label{eqpath2}   
q(P_n) &=
  \sum_{r=0}^{\lfloor n/2\rfloor}
  \left\{ \binom{n-r}{r} + \binom{n-r-1}{r} \right\} x^{r+1} ,
\\ \intertext{and, for $n \ge 0$ and with $y=\sqrt{1+4x}$,}
\label{eqpath3}   
q(P_n)(x) &=
 \frac{(3+y)(y-1)}{4y} \left(\frac{1+y}{2}\right)^{n+1} +
 \frac{(3-y)(y+1)}{4y} \left(\frac{1-y}{2}\right)^{n+1} .
\end{align}
\end{example}

\begin{corollary}
\label{Fib-path}   
For the path $P_n$,
$q(P_n)(1) = F_{n+2}$,
the $(n+2)$'nd Fibonacci number
(with the convention $F_0=0$, $F_1=1$).
\end{corollary}

\begin{example}
\label{cycle}   
The interlace polynomials of the cycles $C_n$ are $q(C_3) = 4x$ and,
for $n \geq 4$ and with $y=\sqrt{1+4x}$,
\begin{align*}
q(C_n)(x) &=
 \left( \frac{1-y}2 \right)^n + \left( \frac{1+y}2 \right)^n
 + \frac{y^4 - 10 y^2 - 7}{16}
 & & \hspace*{-0.4in} \text{for $n$ even,} \\
q(C_n)(x) &=
 \left( \frac{1-y}2 \right)^n + \left( \frac{1+y}2 \right)^n
 + \frac{y^2-5}4
 & & \hspace*{-0.4in} \text{for $n$ odd.} \\
\end{align*}
\end{example}

\section{Partitions, pairings, digraphs, and the polynomial}
\label{sec:facts}   
In this section we touch lightly on how the polynomial partitions
the set of graphs, and related issues.

\subsection{Partition of graphs by the interlace polynomial}

As one would expect,
the interlace polynomial does not resolve graph isomorphism.
A variety of counterexamples is provided by taking two Euler circuits
$\ckt_1$ and $\ckt_2$ in a \tito digraph~$D$,
and producing their corresponding interlace graphs $\HH_1$ and $\HH_2$;
generally $\HH_1$ and $\HH_2$ will not be isomorphic, but
$q(\HH_1) = q(\HH_2)$ since both count circuit partitions of~$D$.
Another set of counterexamples is provided by Remark~\ref{G1G2},
noting that $q(K_m \cup K_{n-m}) = 2^{n-2} x^2$ regardless of $m$.
Another counterexample is provided by the 5-cycle,
and the 5-cycle with one chord-edge added:
both have the polynomial $6x + 5x^2$.

In fact, the polynomial does not always even distinguish trees.
For example, the path 1--2--3--4--5--6--7
with branches 3--8 and 5--9
has polynomial $2x + 9x^2 + 17x^3 + 13x^4 + 4x^5$,
and so does the path 1--2--3--4--5--6
with branches 4--7--8 and 5--9.

We have not seriously explored either question, but
we imagine that the number of different polynomials given by
graphs (respectively, trees) of order $n$ is exponentially large in $n$,
but that the chance that a graph (tree) is uniquely identified by
its polynomial is exponentially small.

\subsection{(Lack of) monotonicity}

Inspection of a few examples initially suggests that
if $G$ and $G'$ are graphs of the same order,
with edge sets $E \subseteq E'$, then $q_G(1) \leq q_{G'}(1)$.
In fact this is \emph{not} the case.
A small example is when $G'$ is the 4-spoke wheel
(with edges consisting of the cycle 1, 2, 3, 4 united
with the star from 5 to 1, 2, 3, and 4),
and $G$ is the same graph with a ``circumference edge'' deleted
(say the edge from 1 to~2).
Here $q_{G}(x) = 6x + 5x^2$
and $q_{G}(x) = 4x + 4x^2 + x^3$,
so $q_{G}(1) = 11$ and $q_{G'}(1) = 9$.

\subsection{Pairings and partitions by digraphs and interlace graphs}
A ``pairing'' may be defined as
a list of vertices $1$ to $n$, in which each vertex appears twice;
we think of this as the order in which a circuit visits these vertices.
From a pairing we can construct a unique \tito digraph,
and also a unique interlace graph.
Since several pairings may lead to the same \tito digraph
or interlace graph, the \tito digraphs and the interlace graphs
both partition the pairings.

Neither of these partitions
\emph{on the pairings}
refines the other.
In a smallest example, the pairings
\mbox{1 1 2 2 3 3} and \mbox{1 1 2 3 3 2}
have isomorphic interlace graphs (both are edgeless),
but non-isomorphic \tito digraphs.
Also in a smallest example, the pairings
\mbox{1 2 3 1 3 4 2 4} and \mbox{1 2 4 1 3 4 2 3}
have non-isomorphic interlace graphs,
but the same \tito digraph.

However, as shown by Propositions \ref{confuse1} and \ref{confuse2},
the interlace graphs and \tito graphs partition each
other: having pairings with the same interlace graph is an equivalence
relation on \tito graphs, and having pairings with the same
\tito graph is an equivalence relation on interlace graphs.

Note first that Lemmas~\ref{ukkpev} and~\ref{lem:transpose}
have the following immediate consequence.
For a graph~$G$, let $\Pset(G)$ be the orbit of $G$ under pivoting.
Then for every $C \in \Cset(D)$, $H(\Cset(D)) = \Pset(H(C))$.

\begin{proposition}
\label{confuse1}   
Let $\Hset(D)$ denote the set of interlace graphs of
all Euler circuits of $D$.
Let $D$ and $D'$ be \tito digraphs such that
$\Hset(D) \cap \Hset(D') \neq \emptyset$.
Then $\Hset(D) = \Hset(D')$.
\end{proposition}

\begin{proof}
Choose $C \in \Cset(D)$ and $C' \in \Cset(D')$ such that
$H(C) = H(C')$.
Then $\Hset(D) = H(\Cset(D)) = \Pset(H(C))
= \Pset(H(C')) = H(\Cset(D')) = \Hset(D')$.
\end{proof}

Let $\Dset(\HH)$ denote the set of \tito digraphs obtained from all
pairings with interlace graph~$\HH$.
Note that for every interlace graph $\HH$ and edge $ab \in E(\HH)$,
$\Dset(\HH) = \Dset(\HH^{ab})$.
(Any $D \in \Dset(\HH)$ is witnessed by a circuit $C$
with $D = D(C)$ and $\HH = \ilace(C)$;
$C^{ab}$ is a witness that $D \in \Dset(\HH^{ab})$.
Thus $\Dset(\HH) \subseteq \Dset(\HH^{ab})$;
by symmetry, then, $\Dset(\HH) = \Dset(\HH^{ab})$.)
Putting it slightly differently, $\Dset(\HH) = \Dset(\Pset(\HH))$.

\begin{proposition}
\label{confuse2}   
If $\HH$ and $\HH'$ are interlace graphs such that
$\Dset(\HH) \cap \Dset(\HH') \neq \emptyset$,
then $\Dset(\HH) = \Dset(\HH')$.
\end{proposition}

\begin{proof}
Let $D \in \Dset(\HH) \cap \Dset(\HH')$ and choose
$C, C' \in \Cset(D)$ such that
$\ilace(C) = \HH$ and $\ilace(C') = \HH'$.
$\Cset(D)$ is a single orbit under transpositions,
and contains $C$ and $C'$,
so $\HH$ and $\HH'$ belong to the same orbit under pivoting:
$\Pset(\HH) = \Pset(\HH')$.
Hence $\Dset(\HH) = \Dset(\Pset(\HH)) = \Dset(\HH')$, as claimed.
\end{proof}

\section{Interlace polynomials of substituted and rotated graphs}
\label{lesssimple}   

We now compute interlace polynomials for
``substituted'' and ``rotated'' graphs.
These are of mild interest in themselves,
and are needed for the proofs (and occasionally the statements)
of the results in Section~\ref{extremal},
on extremal properties of the interlace polynomial.

\subsection{Substituted graphs}
\label{sec:substitution}   

The following definition of graph substitution is recalled from~\cite{BolMGT}.
\begin{definition}
\label{def:solid}   
If $G$ is a graph with vertices $v_1,\ldots,v_n$,
and $G_1,\ldots,G_n$ are arbitrary graphs on disjoint vertex sets,
then $\Gstar = G[G_1,\ldots,G_n]$ is obtained from $\bigcup_{i=1}^n G_i$
by joining all vertices of $G_i$ to all vertices of $G_j$ whenever
$ij \in E(G)$.
\end{definition}

We say that $\Gstar$ is obtained from $G$ by
\emph{substituting} $G_1,\ldots,G_n$ for the vertices
or by \emph{replacing} the vertices with $G_1,\ldots,G_n$.
Note that if we replace the vertices of $G$ one by one
with the graphs $G_1,\ldots,G_n$, we get the same graph $\Gstar$.
If each $G_i$ is a complete graph then we call
$\Gstar$ a \emph{solid graph with template $G$},
or a {\em solid $G$-graph}.
Similarly, if each $G_i$ is an edgeless graph
then $\Gstar$ is a \emph{thick graph with template $G$}
or a {\em thick $G$-graph}.%
\footnote{The case that a substituted graph $G_i$ is
the null graph ($\order{G_i}=0$) is unnatural but sometimes convenient,
for example in the proof of Proposition~\ref{secondmax}.}

Clearly, a graph $G$
is a solid graph with a template of order $k$ if and only if
$V(G)=V_1\cup\ldots\cup V_k$ such that each induced subgraph
$G \restricted{V_i}$ is complete and for every
pair $i, j$, $1\le i < j \le k$, either $G$ contain all $V_i-V_j$
edges, or it contains none of them. Thick graphs have a similar
description, except $G \restricted{V_i}$ is not complete but edgeless.
In particular,
the complement $\overline G$ of a solid graph $G$ with template $H$ is
a thick graph with template $\overline H$. Specializing further, the
complement of a solid path of length $3$ is a thick path of length $3$.

Immediately from the definition, substitution is transitive:
if $\Gstar$ has template $G$ and $G$ has template $H$,
then $\Gstar$ has template $H$.
Specifically, we have the following.

\begin{proposition}
If $\Gstar = G[G_1,\ldots,G_n]$ and $G=H[H_1,\ldots,H_m]$
then there exist graphs $H'_1,\ldots,H'_n$ such that
$\Gstar = H[H'_1, \ldots, H'_m]$.
\end{proposition}

\begin{proof}
The result has no mathematical content, but is confusing.
Each vertex $v_i$ in $H$ is replaced by a graph $H_i$ to form $G$.
Then each vertex in $G$ --- write it as $u_i^j$ to indicate that it
is the $j$'th vertex within graph $H_i$ --- is itself substituted
with a graph $G_i^j$.
Writing $n_i = \order{H_i}$, the net result is, specifically,
\begin{align*}
\Gstar = H[ H_1[G_1^1,\ldots,G_1^{n_1}], \ldots, H_m[G_m^1,\ldots,G_m^{n_m}] ] .
\end{align*}
\end{proof}

In the case of greatest interest, we have the following.

\begin{lemma}
If $\Gstar$ is a solid $G$-graph
and $G$ is a solid $H$-graph,
then $\Gstar$ is a solid $H$-graph.
\end{lemma}

\begin{proof}
This is a specialization of the previous proposition in which
each $H_i$ is a complete graph and each $G_i^j$ is a complete graph,
and thus each $H_i[G_i^1,\ldots,G_i^{n_i}]$ is a complete graph,
making $\Gstar$ a solid $H$-graph.
\end{proof}

\begin{lemma}
\label{lemma:solid}   
Let $a$ be a vertex of a graph $H$ such that $G=H-a$ is a solid graph
with template of order $k$. Then $H$ is a solid graph with template of
order at most $2k+1$.

Also, if $ab$ is an edge of a solid graph $G$ with template of order~$k$,
then $\Gab$ is also a solid graph with template of order $k$.
\end{lemma}

\begin{proof}
First, let $V_1, \ldots V_k$ be the vertex
classes showing that $G=H-a$ is a solid graph with template of order
$k$. For $1\le i \le k$, let $U_i=V_i \cap \Gamma (a)$ and
$W_i=V_i-\Gamma (a)$. Then the classes $\{a\}, U_1, W_1, \ldots , U_k,
W_k$ (some of which may be empty), show that $H$ is a solid graph with
template of order $k$.

Second, if $a$ and $b$ belong to the same class $V_i$ then
$\Gab=G$. If, on the other hand, $a\in V_1$ and $b\in V_2$, say, then
$V_1'=V_1\cup \{b\}-\{a\}$, $V_2'=V_2\cup \{ a \} - \{ b \}$,
$V_3$, \ldots, $V_k$ show that $\Gab$ is a solid graph with template
of order $k$.
\end{proof}

\begin{proposition}
\label{prop:substitution}   
Let $\Gstar = G[K_{m_1},\ldots, K_{m_n}]$.
Then $q(\Gstar) = q(G) \cdot 2^{\order{\Gstar}-\order{G}}$.
\end{proposition}

\begin{proof}
For any $m_i>1$, let $H=G[K_{m_1},\ldots,K_{m_i-1},\ldots, K_{m_n}]$.
Let $a,b$ denote two of the vertices
in the $K_{m_i}$, and pivot-reduce on edge~$ab$.
Note that $(\Gstar)^{ab} = \Gstar$ and that
$\Gstar-a = (\Gstar)^{ab}-b = H$, so that $q(\Gstar) = 2q(H)$.
The proposition follows inductively.
\end{proof}

The simplest case of graph substitution is
duplicating a single vertex $a$,
\ie, replacing $a$ with $E_2$; we will write this $\Gstar = G \circ a$.
For the case of multiplying the $i$'th vertex
$k_i$-fold, \ie, taking $G[E_{k_1}, \ldots, E_{k_n}]$,
we will write $\Gstar=G[k_1,\ldots,k_n]$.
As far as the interlace polynomial is concerned,
general graph substitution can be expressed in terms
of vertex multiplication, via the following theorem.

\begin{proposition}
For the general graph substitution
$\Gstar=G[G_1,\ldots,G_n]$,
with $q(G_i) = \sum a(i,k) x^{k}$, we have
$$
q(\Gstar)
 =
\sum_\csub{k_1,...,k_n} a(1,k_1) \cdots a(n,k_n) \:
q(G[k_1,\ldots,k_n]).
$$
\end{proposition}
\begin{proof}
In $G[G_1,\ldots,G_n]$, the vertices in $G_i$ cannot distinguish among any
vertices in $G_j$, regardless of whether $ij$ is an edge in $G$.
Thus, pivots
can be performed within each $G_i$ separately, with no interactions.
\end{proof}

The interlace polynomials $q(G[k_1,\ldots,k_n])$  occurring above
can all be expressed as combinations of the the $2^n$ interlace
polynomials $q(H)$, as $H$ runs over the induced subgraphs of $G$.
We carry this out in the next
three propositions, first calculating the effect of
a single vertex duplication, then calculating the effect of
duplicating an arbitrary subset of the vertices, and finally
handling the general case.

\begin{proposition}
\label{prop dup}   
If $a$ is a vertex of a graph $G$, then
$$
q(G \circ a) = (1+x) q(G) - xq(G-a).
$$
\end{proposition}
\begin{proof}
If $a$ is an isolated vertex then
$q(G \circ a) = x q(G)$ and $x q(G-a) = q(G)$, proving the assertion.
Otherwise, pivot on an edge $ab$.
Consider the ``duplicate'' vertex $a'$
in the context of the four classes of Definition~\ref{def:pivot}:
neighbors $N(a)$ of $a$ alone,
neighbors $N(b)$ of $b$ alone,
neighbors $N(ab)$ of both,
and vertices $N_0$ that are neighbors of neither.
We know that before pivoting, like $a$:
$a' \in N(b)$,
$a'$ is adjacent to every vertex in $N(a)$ and $N(ab)$,
and $a'$ is adjacent to no vertex $N(b)$ nor $N_0$.
After pivoting, then,
$a' \in N(b)$,
$a'$ is now adjacent to no vertex in $N(a)$ nor in $N(ab)$,
and $a'$ remains adjacent to no vertex in $N(b)$ nor $N_0$.
So in $\Gab$, $a'$ is adjacent only to $b$,
and in $\Gab-b$, $a'$ is an isolated vertex:
$(G \circ a)^{ab} - b = (\Gab - b) \cup E_1$.
Thus
$q(G \circ a)
 = q(G \circ a - a) + q((G \circ a)^{ab} - b)
 = q(G) + x q(\Gab - b)
 = q(G) + x [q(G) - q(G-a)]
 = (1+x) q(G) - x q(G-a)
$.
\end{proof}

\begin{proposition}
\label{prop m dup}   
For a graph $G$
with vertices $v_1,\ldots,v_n$, for any $1 \leq m \leq n$,
the interlace polynomial for the graph $\Gstar$ formed
by duplicating $v_1,\ldots,v_m$ can be expressed in terms of the
interlace polynomials of the  $2^m$ graphs formed by deleting some subset
of $\{ v_1,\ldots,v_m\}$, as follows.
Let $k_1=\cdots =k_m=2$,
$k_{m+1}=\cdots=k_n=1$, and $\ell_{m+1}=\cdots=\ell_n=1$. Then
with $p_1(x)=(1+x)$ and $p_0(x)=-x$,
\begin{align}
\label{m dup}   
q(G[k_1, \ldots, k_n]) &=
\sum_{\ell_1=0}^1 \cdots \sum_{\ell_m = 0}^1
 \left\{
 q(G[\ell_1,\ldots,\ell_n])
\prod_{i=1}^m p_{\ell_i}(x)
 \right\}  .
\end{align}
\end{proposition}

\begin{proof}
This proposition is proved by repeatedly applying
Proposition~\ref{prop dup}.
Note first that for $m=1$, with $v_1=a$, \eqref{m dup} is precisely
Proposition~\ref{prop dup}.
Similarly, for $m=2$, writing $v_1=a$, $v_2=b$, \eqref{m dup} simply says
\begin{equation}
\label{m=2}   
q(G \circ a \circ b) = (1+x)^2 q(G)
 - x(1+x)(q(G-a)+q(G-b))  + x^2 q(G-a-b).
\end{equation}
Let $H=G \circ a$, so that $H \circ b = G \circ a \circ b$
has $q(H \circ b) = (1+x)q(H) - x q(H-b)$.  Note that
$H-b = (G-b) \circ a$ has $q(H-b) =$
$q((G-b) \circ a) = $
$ (1+x) q(G-b) - x q (G-a-b)$.
Thus
\begin{align*}
q(G \circ a \circ b) &= q(H \circ b) \\
 & = (1+x) [q(H)] - x [q(H-b)]   \\
 & = (1+x) [ (1+x) q(G) - x q(G-a) ]  - x [(1+x) q(G-b) - x q (G-a-b)],
\end{align*}
which shows \eqref{m=2}.  The general case
follows by induction on $m$,
applying Proposition \ref{prop dup} to compute the effect of
duplicating the $m+1$'st, in the graph where the first $m$
vertices have already been duplicated.
\end{proof}

\begin{proposition}
\label{prop mult}   
Given a graph $G$ on vertices $v_1,\ldots,v_n$, with $n \geq 1$,
define polynomials $q(k_1,\ldots,k_n)(x)$ by
\begin{align*}
q(k_1,\ldots,k_n)(x)=
\sum_{\ell_1=0}^1 \cdots \sum_{\ell_n = 0}^1 (-1)^{n+\ell}
 q(G[\ell_1,\ldots,\ell_n])(x)
 \prod_{i=1}^n
 \sum_{j=1-\ell_i}^{k_i-1} x^j
\end{align*}
where $\ell = \sum_{i=1}^n \ell_i$.
Then for $k_1,\ldots,k_n \geq 1$,
\begin{equation}
\label{general dup}   
q(G[k_1, \ldots, k_n])(x) = q(k_1,\ldots,k_n)(x) .
\end{equation}
\end{proposition}

\begin{proof}
The proof is by induction on $k_1+\cdots+k_n$, and the basis
consists of $2^n$  cases --- those with all
$k_i \leq 2$. (The case
of \eqref{general dup}
with all $k_i=1$ is trivial, and the other basis
cases are given by Proposition \ref{prop m dup}.)
Assume then that some $k_i \geq 3$; without loss of generality
we may take $i=1$.
Write $a=v_1$ and $H=G[k_1-1,k_2, \ldots, k_n]$, so that by
Proposition \ref{prop dup},
$q(G[k_1,k_2, \ldots, k_n])(x)=q(H \circ a)
= (1+x) q(H) - x q(H-a)$.  By induction, $q(H)$ and
$q(H-a)$ may be expanded as $q(k_1-1,\ldots,k_n)$ and
$q(k_1-2,k_2,\ldots,k_n)$, in which, for each
choice of $\ell_1,\ldots,\ell_n$, all factors are
equal, except for the factor indexed by $i=1$.  For this
factor, the contribution to $(1+x) q(H) - x q(H-a)$ is
$(-1)^{1-\ell_1}$ times
$$
(1+x)
  \left(
 \sum_{j=1-\ell_1}^{(k_1-1)-1} x^j \right)
-
(x)
 \left(
 \sum_{j=1-\ell_1}^{(k_1-2)-1} x^j \right)
\ = \ \
  \left(
 \sum_{j=1-\ell_1}^{k_1-1} x^j \right) ,
$$
which is the coefficient of
$(-1)^{1-\ell_1}$ in the corresponding
factor in the expression for
$q(k_1,k_2,\ldots,k_n)$.  This shows
$q(G[k_1,k_2, \ldots, k_n])(x)=q(k_1,k_2,\ldots,k_n)(x)$,
completing our proof by induction.
\end{proof}

\begin{corollary}
For $r \geq 1$ and $k_1,\ldots,k_r \geq 1$, the complete
$r$-partite graph $\Gstar = K_r(k_1,\ldots,k_r)$ has
\begin{equation}
\label{complete r partite}  
q(\Gstar)(x) =
 \frac{x}{2} \prod_{i=1}^r (2 + x + \cdots + x^{k_i-1})
  + (-1)^r (1-\frac{x}{2}) \prod_{i=1}^r (x + \cdots + x^{k_i-1}) .
\end{equation}
\end{corollary}

\begin{proof}
With $n \equiv r$, this is the special case
 $G=K_n$ in Proposition \ref{prop mult}.  In detail, for
$l_1,\ldots,l_n \in \{ 0,1 \}$ with $s=l_1+\cdots+l_n$
note that
$G(l_1,\ldots,l_n)=K_s$, with $q(K_s)(x)=2^{s-1}x$
\emph{provided that} $s>0$, but
$q(K_s)(x)=1= 2^{s-1}x + (1-x/2)$ for $s=0$.
Adding this ``correction'' to the term with $l_1=\cdots=l_n=0$
gives a sum which factors:
\begin{align*}
\mbox{} & \hspace*{-1cm} q(\Gstar)+
  (1-\frac{x}{2}) \prod_{i=1}^n  \left(  \sum_{j=1}^{k_i-1} x^j \right)
\\ &
= \sum_{\ell_1=0}^1 \cdots \sum_{\ell_n = 0}^1
 \left\{
\frac{x}{2} \
2^{l_1+\cdots+l_n}
 \prod_{i=1}^n \left[ \left(
 \sum_{j=1-\ell_i}^{k_i-1} x^j \right) (-1)^{1-\ell_i} \right]
 \right\}
\\ &
=\frac{x}{2} \
 \prod_{i=1}^n \left[ 2^1 \left(1 + x + \cdots + x^{k_i-1} \right)
- 2^0 \left( x + \cdots + x^{k_i-1} \right) \right]
\\ &
=\frac{x}{2} \prod_{i=1}^r (2 + x + \cdots + x^{k_i-1}).
\end{align*}
\end{proof}

\subsection{Rotated graphs}
\label{rotation}   

Let $H$ be a graph including distinct vertices
$u,v,w$, in which $uw$ is an edge
and
no other edges on $w$ are allowed.
Let $G$ be obtained from $H-w$ by toggling the $uv$
relation, so that $uv$ is an edge of $G$ if and only if it
is not an edge of $H$.
We say that $(G,H)$ is a {\em rotation}.

The operation's name comes from the case where $uv$ was
an edge in $G$, and this edge is ``rotated'' around
the pivot point $u$ to become the new edge $uw$.

\begin{theorem}
\label{rotate-geq}   
Let $(G,H)$ be a rotation and $1 \leq x$.
Then $q_G(x) \leq q_H(x)$.
\end{theorem}

\begin{proof}
We prove the inequality by induction on the order of $H$.
For $|H|=3$ there are
two
cases to check,
according to whether $uv$ is or is not an edge
in $H$.
When $uv$ is an edge, $q_H(x) = 2x+x^2 \geq x^2 = q_G(x)$,
as long as $x \geq 0$.
When $uv$ is not an edge,
$q_H(x) = 2x^2 \geq 2x = q_G(x)$ as long as $x \geq 1$.

In proving the induction step,
fix an $x \geq 1$,
and assume that $|H|>3$.
We consider the following cases. Other than $u$ and $v$, $G$ has:
(1) only isolated vertices;
(2) an edge $ab$, $a$ and $b$ distinct from $u$ and $v$;
(3) no edge as above, and some vertex $c$ adjacent to $v$ but not $u$;
(4) none of the preceding, and some vertex $d$ adjacent to $u$ but not $v$;
(5) none of the preceding, and some vertex $e$ adjacent to both $u$ and $v$.

If there is an isolated vertex $a$ of $G$, then $a$ is also isolated in $H$, and
$q(G;x)=xq(G-a;x) \leq xq(H-a;x) = q(H;x)$,
where the inequality follows from the inductive hypothesis.

If there is an edge $ab$ independent of $u$ and $v$,
pivot-reduce both $G$ and $H$ on~$ab$.
The $ab$-pivot either toggles the edge $uv$ in both $G$ and in $H$,
or leaves it unchanged in both $G$ and $H$;
also it does not affect the neighborhood of $w$.
Thus $(\Gab,\Hab)$ is a rotation,
and so is $(\Gab-b,\Hab-b)$.
Trivially, $(G-a,H-a)$ is also a rotation, and, by the inductive hypothesis,
$q(G;x) = q(G-a;x)+q(\Gab-b;x) \leq q(H-a;x)+q(\Hab-b;x) = q(H;x)$.

Failing the above, it must be that every other vertex in $G$ is
adjacent to $u$ or $v$ or both, but to no other vertices.
Suppose there is a vertex $c$ adjacent only to $v$, so $G^{cv}=G$.
Then $q(G;x) = q(G-c;x)+q(G-v;x)$.
By the inductive hypothesis, $q(G-c;x) \leq q(H-c;x)$.
Also, $q(H-v;x) = q(H-v-w;x) + q(H-v-u;x) \geq q(H-v-w;x) = q(G-v;x)$.
Summing the two inequalities, $q(G;x) \leq q(H;x)$.

The next case is that there is a vertex $d$ adjacent only to $u$.
Then $q(G;x) = q(G-d;x)+q(G-u;x) \leq q(H-d;x)+q(H-u;x) = q(H;x)$:
$q(G-d;x) \leq q(H-d;x)$ is by the inductive hypothesis,
while $q(H-u;x) = xq(H-u-w;x) = xq(G-u;x) \geq q(G-u;x)$ because
$H-u$ and $G-u$ differ just in the isolated vertex~$w$.

The final case is that every vertex in $G$, other than $u$ and $v$,
is adjacent to both $u$ and~$v$.
Pick one such vertex $e$ and pivot-reduce $G$ and $H$ on $ev$.
By the inductive hypothesis, $q(G-e;x) \leq q(H-e;x)$.
On the other hand, $H^{ev}$ consists of a star on $v$,
an edge $uw$, and possibly an edge $uv$; $G^{ev}$ is the same
but without the vertex $w$.
Pivot-reducing $H^{ev}$ on $wu$ gives
$q(H^{ev}) = q(H^{ev}-w) + q(H^{ev}-u) = q(G^{ev}) + xq(G^{ev}-u) \geq q(G^{ev})$.
Summing the two inequalities, here too
$q(G;x) \leq q(H;x)$, and we are done.
\end{proof}

\section{Extremal properties of the interlace polynomial}
\label{extremal}   

In this section we
consider extremal values of the interlace polynomial's degree,
its evaluation at~1, and its number of nonzero terms,
in terms a graph's order $|G|$ (the number of vertices)
and its size $e(G)$ (the number of edges).

\subsection{Extremal values of the degree of $\mathbf{q(G)}$}

\begin{proposition}
\label{maxdeg}   
For any graph $G$ of order $n$, $\deg(q(G)) \leq n$,
with equality only for $G=E_n$.
Also, $q_G(1) \geq 1$, with equality only for $G=E_n$.
\end{proposition}

\begin{proof}
The first assertion follows from Remark~\ref{indep}: $\deg(q(G))=n$
means the independence number $\alpha(G)=n$, and $G$ must be an edgeless graph.

The second assertion is proved by induction on the order of $G$.
If $G=E_n$, then $q_G(1)=1$ and we are done.
Otherwise there is an edge $ab$ in $G$,
and $q(G) = q(G-a) + q(\Gab-b) \geq 2$ by the inductive hypothesis.
\end{proof}

\begin{remark}
If $q(G)$ is purely linear, then $G$ is a complete graph.
\end{remark}

\begin{proof}
By Remark~\ref{indep}, $\deg(q(G))=1$ means $\alpha(G)=1$,
and $G$ is complete.
\end{proof}

\subsection{Extremal values of $\mathbf{q_G(1)}$}

Propositions
\ref{sizeLB}, \ref{sizeUB},
\ref{orderLB}, and \ref{orderUB}
provide lower and upper
bounds for $q_G(1)$ in terms of $G$'s size and order. 
The proof of Proposition~\ref{sizeLB} will need the following lemma.

\begin{lemma}
\label{Fib-tree}   
For any tree $G$ of order $n$, $q_G(1) \leq F_{n+1}$,
with equality achieved precisely by a path $P_{n-1}$.
\end{lemma}

\begin{proof}
One direction is immediate from Corollary~\ref{Fib-path}:
if $G = P_{n-1}$ then $q_G(1) = F_{n+1}$.
For the other direction we apply induction on~$n$
to show that any tree $G$ achieving the bound must be a path.
The base case with $n=1$ is trivial,
as is that $G \neq E_n$.

If $n>1$, then pick any leaf $a$, with parent $b$.
The pivot reduction on $ab$ gives
$q(G) = q(G-a) + q(G-b) = q(G-a) + xq(G-b-a)$.
Evaluating at $x=1$,
$q(G;1) = q(G-a;1) + q(G-b-a;1) \leq F_n + F_{n-1}$
by the inductive hypothesis, with equality only if $G-a$ and
$G-a-b$ are both paths.  For this to be the case, $b$ must be
a terminal of the path $G-a$, and therefore $G$ is also a path.
Hence the only tree achieving the upper bound is the path.
\end{proof}

For convenience, we subsume complete bipartite graphs
and edgeless graphs into complete tripartite graphs,
by allowing one or two of the three vertex classes to be empty.
\begin{proposition}[Size, lower bound]
\label{sizeLB}   
For every graph $G$, $q_G(1) \geq \size{G}+1$,
with equality if and only if $G$ consists of a
complete tripartite graph and isolated vertices.
\end{proposition}

\begin{proof}
Apply induction on the order $n$ of $G$.
For $n \leq 3$ or $\size{G}=0$ the result holds by inspection,
so assume that $n \geq 4$, $\size{G}>0$, and the assertion holds for all smaller
values of $n$.

Let $ab \in E(G)$.  Then
\begin{align}
q(G;1) &= q(G-a;1) + q(\Gab-b;1)        \notag
\\ & \geq [(\size{G}-d(a))+1] + [\size{\Gab-b}+1] \notag
\\ & \geq [(\size{G}-d(a))+1] + [(d(a)-1)+1]  \notag
\\ &= \size{G}+1 , \notag
\end{align}
where the the first inequality
follows from the inductive hypothesis.
This proves the inequality and ---
since both inequalities hold with equality when $G$ consists of a
tripartite graph and isolated vertices ---
also proves the ``if'' direction for equality.

If equality holds, then $\size{\Gab-b} = d(a)-1$,
thus all edges of $\Gab-b$ are incident to $a$,
so all edges of $\Gab$ are incident to $a$ or $b$.
Write $\Gamma$ for the neighborhood in $\Gab$.
Letting $V_a = \Gamma(a) \backslash \Gamma(b)$,
$V_b = \Gamma(b) \backslash \Gamma(a)$,
and $V_{ab} = \Gamma(a) \cap \Gamma(b)$,
pivoting $\Gab$ on $ab$ to recover $G$,
we find that $G$ is the complete tripartite graph with vertex
classes $V_a$, $V_b$, and $V_{ab}$.
(Note that $V_a$ and $V_b$ are both non-empty, since
$b \in V_a$ and $a \in V_b$.)
\end{proof}

The next extremal proposition,
bounding $q_G(1)$ by a Fibonacci number,
uses a property of graph rotation from Section~\ref{rotation}.

\begin{proposition}[size, upper bound]
\label{sizeUB}   
For any connected graph $G$ with $m$ edges, $q_G(1) \leq F_{m+2}$,
with equality achieved precisely by a path $P_m$.
If $G$ consists of components $G_i$ with sizes $m_i$,
then $q_G(1) \leq \prod F_{{m_i}+2}$.
Always, if $G$ has $m$ edges then $q_G(1) \leq 2^m$,
with equality if and only if $G$ consists of $m$ independent edges,
and isolated vertices.
\end{proposition}

\begin{proof}
The proposition's essence is its first assertion.
One direction is immediate from Corollary~\ref{Fib-path}:
if $G = P_m$ then $q_G(1) = F_{m+2}$.

For the other direction, suppose that $q_G(1) \geq F_{m+2}$.
If $G$ is a tree,
by Lemma~\ref{Fib-tree},
it can only achieve the bound
$q(G;q) = F_{n+1} = F_{m+2}$ if it is a path.
Otherwise, by the explicit calculations in Example~\ref{cycle},
$G$ cannot be a simple cycle,
and so it has some vertex $v$ of degree at least~3.
Choose any edge $ab$ in a cycle of $G$ and not incident to $v$,
and ``rotate'' the edge out;
that is, delete it, and add a new vertex incident to $a$ alone.
Repeat until no cycles remain.
As we only deleted cycle edges, the graph remains connected.
By Theorem~\ref{rotate-geq}, each rotation can only increase $q(1)$,
yet the final tree, whose vertex $v$ shows it not to be a path,
has $q(1) < F_{m+2}$.  Thus $q_G(1) < F_{m+2}$,
completing the proof of the main assertion.

The result for graphs with several components follows immediately
by Remark~\ref{G1G2}.
Also, $F_{m+2} \leq 2^m$, with equality iff $m=0$ or~$1$,
proving the final statement.
\end{proof}

\begin{proposition}[order, lower bound]
\label{orderLB}   
Let $G$ be a graph of order $n$, having no isolated vertices.
Then $q_G(1) \geq n$, with equality if and only if $G$ is a star,
or $n=4$ and $G$ consists of two independent edges.
\end{proposition}

\begin{proof}
The ``if'' direction is trivial.
To prove the ``only if'', we apply induction on~$n$.
For $n \leq 4$ the result is easily checked,
so let us assume that $n>4$ and the result holds for
all smaller values of $n$.
If $G$ is disconnected, say $G = G_1 \cup G_2$ with
$\order{G_1} = n_1 \geq 2$ and $\order{G_2} = n_2 \geq 2$,
then, by the inductive hypothesis,
$q_G(1) \geq n_1 n_2 > n$.
If on the other hand $G$ is connected,
then either $G$ is a tree with $n-1$ edges, or $\size{G} \geq n$.
Hence, by Proposition~\ref{sizeLB}, $q_G(1) \geq n$, with equality
iff $G$ is a tree and is an extremal graph of that Theorem.
A graph that is tree and is a complete tripartite graph is a star,
completing the proof.
\end{proof}

This simple proposition was used by
Balister, Bollob\'as, Riordan and Scott \cite{BBRS}
to give a particularly simple proof
of the classical theorem of Bankwitz \cite{Bankwitz}
that if a knot has a nontrivial reduced alternating
diagram then it is nontrivial.

Note that the final proposition in this set of four
has a nice parallel to Lemma~\ref{Fib-tree}.

\begin{proposition}[order, upper bound]
\label{orderUB}   
For any graph $G$ of order $n$,
$q_G(1) \leq 2^{n-1}$, with equality if and only if $G=K_n$.
\end{proposition}

\begin{proof}
The proof is by induction on $n$, with trivial base case $n=1$.
For $n>1$, $G$ cannot be the edgeless graph, so choose an edge $ab$ in $G$.
Then $2^n = q(G) = q(G-a) + q(\Gab-b) \leq 2^{n-1} + 2^{n-1}$,
so both $q(G-a)$ and $q(\Gab-b)$ must achieve their extremal values,
and by the inductive hypothesis each must be isomorphic $K_{n-1}$.
Since $\Gab-b$ is complete, the neighborhood of $a$ in $\Gab$
consists of all vertices ($a$ is connected to $b$ by definition).
The neighborhoods of $a$ in $G$ and in $\Gab$ are always the same
by definition of the pivot operator, thus $a$ is connected to every
other vertex in~$G$.
Since $G-a \isom K_{n-1}$, $G \isom K_n$.
\end{proof}

\subsection{The several largest values of $\mathbf{q_G(1)}$}

Proposition~\ref{orderUB} allows an extension, given
by Proposition~\ref{secondmax}.
It and
the following results rely on properties of graph substitution
from Section~\ref{sec:substitution}.

\begin{proposition}
\label{secondmax}   
Over graphs $G$ of order $n$,
the second-maximum value of $q_G(1)$ is $\frac34 2^{n-1}$,
and is achieved exactly when $G$ is a solid $P_2$-graph
but not a complete graph.
(That is, $G$ is composed of a pair of complete graphs
sharing some but not all their vertices.)
\end{proposition}

\begin{proof}
Observing that a complete graph is a special case of a solid $P_2$-graph,
and is covered by Proposition~\ref{orderUB},
the present proposition is equivalent to:
$q_G(1) \geq \frac38 2^n$ precisely when $G$ is a solid
$P_2$-graph (where one of the three components may be of size~0).

With a little work, we can check that the assertion holds for $n \leq 7$,
so assume that $n > 7$ and the assertion holds for graphs of orders up to $n-1$.
Let $G$ be a graph of order $n>7$ with $q_G(1) \geq \frac38 2^n$.
Then $G$ is not edgeless; let $ab\in E(G)$.
Since $q(G)(1)= q(G-a)(1)+q(\Gab-b)(1)$,
either $q(G-a)(1) \geq \frac38 2^{n-1}$
or $q(\Gab-b)(1) \geq \frac38 2^{n-1}$.
But then, by the inductive hypothesis,
either $G-a$ or $\Gab-b$ is a solid $P_2$-graph.
By the first part of Lemma~\ref{lemma:solid},
either $G$ or $\Gab$ is a solid graph of order $n$
with a template $H$ of order $h \leq 7$.
By Proposition~\ref{prop:substitution},
we have $\frac38 2^n \leq q(H)(1) \cdot 2^{(n-1)-h}$,
so $q(H)(1) \geq \frac38 2^h$
and thus $H$ is a solid $P_2$-graph.
But then so is at least one of $G$ and $\Gab$.
In the second case, $G$ too is a solid $P_2$-graph,
by the second part of Lemma~\ref{lemma:solid}.
(If the template of $G$ were any order-3 graph other than $P_2$,
the template would still have to be connected --- since $\Gab$ is
and therefore $G$ is --- so it could only be $K_3$, which as a
template may be regarded as a special case of $P_2$
two of whose three groups have size~0.)
Hence $G$ is a solid $P_2$-graph, as claimed.
\end{proof}

In fact, the proof above gives the following assertion.
Suppose that every
graph of order $n \le 2k+1$ with $q_G(1) \ge c 2^n$ is a solid
$G_i$-graph for some $i= 1, \ldots, \ell$, where $G_1, \ldots ,
G_\ell$ are graphs of order at most $k$.
Let now $G$ be a graph of order $n$ with $q_G(1)\ge c
2^n$. Then $G$ is a solid $G_i$-graph for some $i, 1\le i \le \ell$.

Just as there is a large gap between $q_G(1)$'s maximum value of $2^{n-1}$
and its second-maximum value of $\frac34 2^{n-1}$,
there seem to be other such gaps.
(Table~6 in~\cite{ABCS00} exhibits some, albeit only for interlace graphs.)
The assertion above might be used to prove the following conjecture.

\begin{conjecture}
There are constants $c_1=\frac12 > c_2 > \ldots$ such that for
every $k \ge 1$ and $n$ sufficiently large, then,
first, there are graphs $G_1, \ldots , G_k$ of order $n$
such that $q(G_i)(1)=c_i 2^n$,
and, second,
if $G$ is any graph of order $n$ with $q_G(1) \geq c_k 2^n$,
then $q_G(1)=c_i 2^n$ for some $i \leq k$.
\end{conjecture}

\subsection{Zero and nonzero terms of $\mathbf{q(G)}$}

\begin{lemma}
\label{solidPath}   
Let $G$ be a connected graph of order at least $\ell+4$ such that every
connected induced subgraph of it with $\ell+3$ vertices is a solid
path of length at most $\ell$.  Then $G$ is a solid path of length at
most $\ell$.
\end{lemma}

\begin{proof}
We may clearly assume that $\ell \ge 3$, for some two vertices $a$
and $b$ the graph $G-a-b$ is a  solid $\ell$-path, and the graphs
$G-a$, $G-b$ and $G-a-b$ are connected. Let $P$ be an induced $\ell$-path in
$G-a-b$. Since the induced subgraph of $G$ with vertex set
$V(P)\cup \{a\} \cup \{ b\}$
is a solid $\ell$-path, so is the entire graph $G$.
\end{proof}

\begin{proposition}
Let $G$ be a graph such that $q(G)$ has precisely two non-zero
terms. Then one of the components of $G$ is a solid path of length $2$
or $3$, and all other components are complete graphs.
\end{proposition}

\begin{proof}
It is easily seen that in proving this we may assume that $G$ is an
incomplete connected graph of order $n\ge 6$, and the assertion holds
for graphs of smaller order. Then every connected incomplete induced
subgraph $H$ of $G$ has precisely two non-zero terms, so  $H$ is a
solid path of length at most $3$.  Writing $\ell$ for the maximal
length of an induced path in $G$, Lemma~\ref{solidPath} implies that
$G$ itself is a solid path of length $\ell$.
\end{proof}

\begin{proposition}
If $\order{G} = n \geq 3$ and $q(G)$ has $n-1$ nonzero terms,
then $q(G) = 2x + x^2 + \cdots x^{n-1}$,
and $G$ is an $n$-star.
\end{proposition}

\begin{proof}
Since by Remark~\ref{order}
$2^n = q_G(2) = \sum_{i=1}^n a_i 2^i$ with each $a_i$ a
non-negative integer, we can only have
$q_G(x) = 2x + x^2 + \cdots x^{n-1}$.
By Remark~\ref{component}, since $q(G)$ has nonzero linear term,
$G$ is connected.
Choose an edge $ab$ of $G$ such that $G-a$ is connected.
In particular, for $n \geq 3$, $G-a$ has some edge, 
so (by Proposition~\ref{maxdeg}) $\deg(q(G-a)) \leq n-2$.
Thus the $x^{n-1}$ term of $q(G)$ must belong to $q(\Gab-b)$.
As $q(\Gab-b)(2) = 2^{n-1}$,
we can only have $q(\Gab-b) = x^{n-1}$,
and so (again by Proposition~\ref{maxdeg}) $\Gab-b = E_{n-1}$.
Since $\Gab$ is connected (see Remark~\ref{pivot-connected}),
it must be the star with center~$b$.
We conclude that $G = G^{(ab)(ab)}$ is also the star with center~$b$.
\end{proof}

\section{Open problems}
\label{sec:open}   

The foregoing was a snapshot of our knowledge of the interlace polynomial
at one time.  As pointed out in the introduction,
the polynomial's understanding can be furthered by its interpretation
as a special case of the Martin polynomial of an isotropic
system \cite{Bouchet-personal-comm00,AignerHolst02},
a topic extensively studied by Bouchet
\cite{Bouchet87-isotropic-systems,
Bouchet88-graphic-presentations,
Bouchet91,
Bouchet-Tutte,
Bouchet-mm3}.
Also, a two-variable generalization of the polynomial
and a clarified view of the one-variable polynomial,
without reference to isotropic systems, are given in~\cite{2poly}.

It is clear that the interlace polynomial is not a special
case of the Tutte polynomial, for one thing because the
interlace polynomial varies over trees of order~$n$,
where the Tutte polynomial is always~$x^{n-1}$.
In \cite{2poly} we show that
a 2-variable extension of the interlace polynomial
has an expansion formula structurally identical to the Tutte polynomial's,
but with a sum taken over vertex subsets rather than edge subsets,
and with ranks computed over $\mathbb{F}_2$ rather than $\mathbb{R}$.
That is, the two polynomials have a close structural relationship
but are completely different objects.

There remain many unanswered questions about the interlace polynomial.
First, we mention a pair of related conjectures.

For the first conjecture, we have numerical evidence
compiled from thousands of random graphs of orders up to~13,
as well as all graphs of orders up to~8.
To interpret its second assertion,
recall from Theorem~\ref{thm:q&r} that when $G$ is an interlace graph
with an associated digraph~$D$, $xq(G;1+x) = r(D;x)$ is the
circuit counting polynomial of~$D$.

\begin{conjecture}
\label{unimodal}   
For any graph~$G$,
the sequence of coefficients of $q_G(x)$
is unimodal:
non-decreasing up to some point, then non-increasing.
The associated polynomial $xq_G(1+x)$
likewise has unimodal coefficients.
\end{conjecture}

This conjecture is reminiscent of an outstanding pair of conjectures
about the chromatic polynomial, namely that its coefficient
sequence is always unimodal \cite{Read68} and in fact
log-concave~\cite{Hoggar74} (see~\cite{Welsh93}).
There are two reasons for caution.
First, a similar conjecture for the Tutte polynomial was
falsified by counterexample in 1993~\cite{Schwarzler93}.

Second, the interlace polynomial's coefficients are certainly not
log-concave; the smallest counterexample is $q(K_{1,3}) = 2x+x^2+x^3$.
(We know of no case where the coefficients of the
associated polynomial $r(x)$ fail to be log-concave.)
This means that,
unlike all our other statements about the interlace polynomial,
it does not suffice to prove Conjecture~\ref{unimodal} for connected
graphs: while a product of log-concave polynomials is log-concave,
a product of unimodal polynomials is not necessarily unimodal.
(Even squaring a unimodal polynomial can result in a polynomial with
an arbitrarily large number of modes, as shown in~\cite{Sato92}.)

We have made little headway in proving the conjecture;
for instance we have not even shown that coefficient sequences
have no internal zeros, \ie, that
the nonzero coefficients form a single block.
And, let us note that the second part of the conjecture is not a
trivial consequence of the first:
it is easy to construct a unimodal polynomial $p(x)$
for which $p(1+x)$ is not unimodal,
for example
$p(x)= 1000 x^3+ x^4 + x^5 +x^6 +x^7 +x^8 + x^9 +x^{10} +x^{11} + x^{12}$.

If $G$ is not the interlace graph of any \tito digraph,
then
$x q_G(1+x)$
has no interpretation as a circuit partition polynomial.
However, Conjecture~\ref{unimodal} would imply that for any
\tito digraph~$D$, the number of partitions $r_k$ into $k$ directed circuits
is unimodal in~$k$.
This leads us to speculate the following.

\begin{conjecture}
\label{circuitsUnimodal}   
For any weakly Eulerian graph or digraph,
the number of partitions into $k$ circuits is unimodal in~$k$.
\end{conjecture}

Beyond Conjectures~\ref{unimodal} and~\ref{circuitsUnimodal} specifically,
there are general questions about circuit counting.
For Eulerian graphs and digraphs,
to what degree are the interesting algebraic properties of the
interlace graph mirrored in such properties of the
circuit partition and Martin polynomials?
(For example, Ellis-Monaghan~\cite{E1}
shows that the Martin polynomial is a translation of
a universal skein-type graph polynomial which is a Hopf map.)
For $k>1$, is there a polynomial-time method for counting
$k$-circuit partitions of graphs, directed graphs,
or even \tito directed graphs?

There are also basic, unanswered questions about the interlace polynomial.
Can $q(G)$ be generalized from graphs to matroids?
Can it be computed efficiently?

An exponential-time computation of $q(G)$ is immediate from the
definition, and we suspect that no polynomial-time computation is
possible.
We have explored the issue only very briefly, but we remark
that it may not be immediate to derive
the hardness of the interlace polynomial from
its connection with the Martin polynomial (Theorem~\ref{thm:q&r})
or the further connection to the diagonal Tutte polynomial,
because the connections apply only for limited classes of graphs
where the latter polynomial may not be hard.
At least some evaluations are easy:
For interlace graphs~$H$, $q_H(1)$ can be
computed in polynomial time:
find a circle arrangement $\ckt$ for $H$~\cite{Spinrad}
to obtain a \tito graph $D$,
and count its Euler circuits with
a combination of the directed matrix-tree
theorem~\cite{KirchoffMT,Tutte48}
(see also~\cite[p.~58]{BolMGT}) and
the so-called BEST theorem~\cite{ST41,BE51} (see
also~\cite[p.~18]{BolMGT}).
By the same token, if the full circuit-counting polynomial $r(D)$ for a
\tito digraph $D$ cannot be computed in polynomial time,
then $q(G)$ cannot be computed in polynomial time.

How many different polynomials $q_G$ are there,
say for graphs of order~$n$?
How many graphs may share a single polynomial?
Is there a way to recognize a polynomial obtainable as $q(G)$?
What can we say about $q(G)$ when $G$ is a random graph?

We think that the real question is what $q(G)$ computes
about $G$ itself (say when $G$ is not an interlace graph).
One referee suggests that ``if $S=(L,V)$ is the isotropic
system defined by a graphic presentation using $G$, say $(G,A,B)$,
then $q(G;x)$ is equal to the
restricted Tutte-Martin polynomial $m(G,A+B;x)$''
and $m(G,A+B,1)$ 
``counts the number $n$ of Eulerian vectors of the
isotropic system $S$ that are supplementary of the vector $A+B$.''
While this is surely an answer to the question, we continue to
hope for more elementary combinatorial interpretations.
(For example in \cite{2poly} it is shown that a specialization
of the two-variable interlace polynomial is the independence
polynomial, but that specialization is not the $q(G)$ considered here.)

\section*{Acknowledgments}
We are grateful to Bill Jackson for bringing the work of
Andr\'e Bouchet to our attention,
and to Bouchet and the anonymous referees for their elucidation
of this polynomial's debt to the Tutte-Martin polynomial,
and many helpful detailed comments.
Thanks also to Hein van der Holst for discussing \cite{AignerHolst02}
with us.

\bibliography{poly}

\providecommand{\bysame}{\leavevmode\hbox to3em{\hrulefill}\thinspace}
\providecommand{\MR}{\relax\ifhmode\unskip\space\fi MR }
% \MRhref is called by the amsart/book/proc definition of \MR.
\providecommand{\MRhref}[2]{%
  \href{http://www.ams.org/mathscinet-getitem?mr=#1}{#2}
}
\providecommand{\href}[2]{#2}
\begin{thebibliography}{dBvAE51}

\bibitem[ABCS00]{ABCS00}
Richard Arratia, B{\'e}la Bollob{\'a}s, Don Coppersmith, and Gregory~B. Sorkin,
  \emph{Euler circuits and {DNA} sequencing by hybridization}, Discrete Appl.
  Math. \textbf{104} (2000), no.~1-3, 63--96. \MR{2001i:92032}

\bibitem[ABS00]{ABS00}
Richard Arratia, B\'{e}la Bollob\'{a}s, and Gregory~B. Sorkin, \emph{The
  interlace polynomial: A new graph polynomial}, Proceedings of the Eleventh
  Annual {ACM--SIAM} Symposium on Discrete Algorithms (San Francisco, CA),
  January 2000, pp.~237--245.

\bibitem[ABS01]{2poly}
\bysame, \emph{A two-variable interlace polynomial},
  http:\-//arxiv.org\-/ps\-/math.CO\-/0209054; submitted for journal
  publication, 2001.

\bibitem[AvdH04]{AignerHolst02}
Martin Aigner and Hein van~der Holst, \emph{Interlace polynomials}, Linear
  Algebra Appl. \textbf{377} (2004), 11--30. \MR{2 021 600}

\bibitem[Ban30]{Bankwitz}
C.~Bankwitz, \emph{\"uber die torsionzahlen der alternierenden knoten}, Math.
  Annalen \textbf{103} (1930), 145--161.

\bibitem[BBCP02]{BBCP00}
P.~N. Balister, B.~Bollob{\'a}s, J.~Cutler, and L.~Pebody, \emph{The interlace
  polynomial of graphs at {$\mathit{-1}$}}, European J. Combin. \textbf{23}
  (2002), no.~7, 761--767. \MR{2003h:05126}

\bibitem[BBD97]{Bouchet-Tutte}
Dominique B{\'e}nard, Andr{\'e} Bouchet, and Alain Duchamp, \emph{On the
  {M}artin and {T}utte polynomials}, Tech. report, D{\'e}partement
  d'Informatique, Universit{\'e} du Maine, Le Mans, France, May 1997.

\bibitem[BBRS01]{BBRS}
P.~N. Balister, B.~Bollob{\'a}s, O.~M. Riordan, and A.~D. Scott,
  \emph{Alternating knot diagrams, {E}uler circuits and the interlace
  polynomial}, European J. Combin. \textbf{22} (2001), no.~1, 1--4.
  \MR{2002c:05055}

\bibitem[Bol98]{BolMGT}
B\'{e}la Bollob\'{a}s, \emph{Modern graph theory}, Graduate Texts in
  Mathematics, vol. 184, Springer, New York, 1998.

\bibitem[Bol02]{Boll02}
B{\'e}la Bollob{\'a}s, \emph{Evaluations of the circuit partition polynomial},
  J. Combin. Theory Ser. B \textbf{85} (2002), no.~2, 261--268.
  \MR{2003f:05066}

\bibitem[Bou87]{Bouchet87-isotropic-systems}
Andr{\'e} Bouchet, \emph{Isotropic systems}, European J. Combin. \textbf{8}
  (1987), no.~3, 231--244. \MR{89b:05066}

\bibitem[Bou88]{Bouchet88-graphic-presentations}
\bysame, \emph{Graphic presentations of isotropic systems}, J. Combin. Theory
  Ser. B \textbf{45} (1988), no.~1, 58--76. \MR{89f:05150}

\bibitem[Bou91]{Bouchet91}
\bysame, \emph{{T}utte-{M}artin polynomials and orienting vectors of isotropic
  systems}, Graphs Combin. \textbf{7} (1991), no.~3, 235--252. \MR{92i:05071}

\bibitem[Bou99]{Bouchet-mm3}
Andr{\'e} Bouchet, \emph{Multimatroids {I}{I}{I}.\ {T}ightness and fundamental
  graphs}, Tech. report, D{\'e}partement d'Informatique, Universit{\'e} du
  Maine, Le Mans, France, September 1999.

\bibitem[Bou00]{Bouchet-personal-comm00}
Andr{\'e} Bouchet, personal communication, September 2000.

\bibitem[dBvAE51]{BE51}
N.G. de~Bruijn and T.~van Aardenne-Ehrenfest, \emph{Circuits and trees in
  oriented graphs}, Simon Stevin \textbf{28} (1951), 203--217.

\bibitem[dF84]{deFraysseix}
Hubert de~Fraysseix, \emph{A characterization of circle graphs}, European
  Journal of Combinatorics \textbf{5} (1984), 223--238.

\bibitem[EM98]{E1}
Joanna~A. Ellis-Monaghan, \emph{New results for the {M}artin polynomial}, J.
  Combin. Theory Ser. B \textbf{74} (1998), no.~2, 326--352. \MR{2000a:05109}

\bibitem[EM99]{E99}
\bysame, \emph{Martin polynomial miscellanea}, Proceedings of the Thirtieth
  Southeastern International Conference on Combinatorics, Graph Theory, and
  Computing (Boca Raton, FL, 1999), vol. 137, 1999, pp.~19--31.
  \MR{2001a:05092}

\bibitem[Gol80]{Golumbic}
Martin~C. Golumbic, \emph{Algorithmic graph theory and perfect graphs},
  Academic Press, New York, 1980.

\bibitem[GSH89]{Gabor}
Csaba~P. Gabor, Kenneth~J. Supowit, and Wen-Lian Hsu, \emph{Recognizing circle
  graphs in polynomial time}, Journal of the ACM \textbf{36} (1989), 435--473.

\bibitem[Hog74]{Hoggar74}
S.~Hoggar, \emph{Chromatic polynomials and logarithmic concavity}, Journal of
  Combinatorial Theory (B) \textbf{16} (1974), 248--254.

\bibitem[Jae88]{J1}
Fran{\c{c}}ois Jaeger, \emph{On {T}utte polynomials and cycles of plane
  graphs}, J. Combin. Theory Ser. B \textbf{44} (1988), no.~2, 127--146.
  \MR{89b:05086}

\bibitem[Kau87]{kauffman87}
Louis~H. Kauffman, \emph{State models and the {J}ones polynomial}, Topology
  \textbf{26} (1987), no.~3, 395--407. \MR{88f:57006}

\bibitem[Kau88]{kauffman88}
\bysame, \emph{New invariants in the theory of knots}, Amer. Math. Monthly
  \textbf{95} (1988), no.~3, 195--242. \MR{89d:57005}

\bibitem[Kir47]{KirchoffMT}
G.~Kirchoff, \emph{{\"Uber die Aufl\"osung der Gleichungen, auf welche man bei
  der Untersuchung der linearen Verteilung galvanischer Str\"ome gef\"uhrt
  wird}}, Ann. Phys. Chem. \textbf{72} (1847), 497--508.

\bibitem[Kot68]{Kot68}
A.~Kotzig, \emph{Eulerian lines in finite $4$-valent graphs and their
  transformations.}, Theory of Graphs (Proc. Colloq., Tihany, 1966), Academic
  Press, New York, 1968, pp.~219--230. \MR{40 \#1298}

\bibitem[LV81]{LV81}
M.~Las~Vergnas, \emph{{E}ulerian circuits of 4-valent graphs imbedded in
  surfaces}, Algebraic Methods in Graph Theory (Szeged, Hungary, 1978), Coll.
  Math. Soc. J. Bolyai, vol.~25, North-Holland, Amsterdam, 1981, pp.~45--477.

\bibitem[LV83]{LV83}
\bysame, \emph{Le polyn\^ome de {M}artin d'un graphe {E}ulerian}, Annals of
  Discrete Math. \textbf{17} (1983), 397--411.

\bibitem[LV88]{LV88}
\bysame, \emph{On the evaluation at (3,3) of the {T}utte polynomial of a
  graph}, Journal of Combinatorial Theory, Series B \textbf{44} (1988),
  367--372.

\bibitem[Mar77]{Martin77}
P.~Martin, \emph{Enum\'erations eul\'eriennes dans les multigraphes et
  invariants de {T}utte-{G}rothendieck}, Ph.D. thesis, Grenoble, 1977.

\bibitem[Mar78a]{M2}
Pierre Martin, \emph{Remarkable valuation of the dichromatic polynomial of
  planar multigraphs}, J. Combin. Theory Ser. B \textbf{24} (1978), no.~3,
  318--324. \MR{58 \#335}

\bibitem[Mar78b]{martin78}
\bysame, \emph{Remarkable valuation of the dichromatic polynomial of planar
  multigraphs}, J. Combin. Theory Ser. B \textbf{24} (1978), no.~3, 318--324.
  \MR{58 \#335}

\bibitem[Pev95]{Pevzner}
Pavel Pevzner, \emph{{DNA} physical mapping and alternating {Eulerian} circuits
  in colored graphs}, Algorithmica \textbf{13} (1995), 77--105.

\bibitem[PTW83]{PolyaNotes}
George P{\'o}lya, Robert~E. Tarjan, and Donald~R. Woods, \emph{Notes on
  introductory combinatorics}, Birkh\"auser Boston Inc., Boston, Mass., 1983.
  \MR{85k:05001}

\bibitem[Rea68]{Read68}
R.C. Read, \emph{An introduction to chromatic polynomials}, Journal of
  Combinatorial Theory \textbf{4} (1968), 52--71.

\bibitem[RR78a]{ReRo76}
R.C. Read and P.~Rosenstiehl, \emph{On the {Gauss} crossing problem},
  Combinatorics, Vol.~2 (Keszthely, Hungary) (A.~Hajnal and V.T. S\'os, eds.),
  Colloquia Mathematica Sociatatis J\'anos Bolyai~18, North Holland, Amsterdam,
  1978, pp.~843--876.

\bibitem[RR78b]{RR}
P.~Rosensthiehl and R.~C. Read, \emph{On the principal edge tripartition of a
  graph}, Advances in Graph Theory (Cambridge Combinatorial Conference, Trinity
  College, Cambridge, 1977) (B.~Bollob\'as, ed.), Ann. Discrete Math., vol.~3,
  1978, pp.~195--226.

\bibitem[Sat92]{Sato92}
Ken-Ito Sato, \emph{Convolution of unimodal distributions can produce any
  number of modes}, Annals of Probability \textbf{21} (1992), no.~3,
  1543--1549.

\bibitem[Sch93]{Schwarzler93}
Werner Schw{\"a}rzler, \emph{The coefficients of the {T}utte polynomial are not
  unimodal}, J. Combin. Theory Ser. B \textbf{58} (1993), no.~2, 240--242.
  \MR{94e:05114}

\bibitem[Spi94]{Spinrad}
Jeremy Spinrad, \emph{Recognition of circle graphs}, Journal of Algorithms
  \textbf{16} (1994), 264--282.

\bibitem[ST41]{ST41}
C.A.B. Smith and W.T. Tutte, \emph{On unicursal paths in a network of degree
  4}, American Mathematical Monthly \textbf{48} (1941), 233--237.

\bibitem[Tut48]{Tutte48}
W.T. Tutte, \emph{The dissection of equilateral triangles into equilateral
  triangles}, Proc. Cambridge Philos. Soc \textbf{44} (1948), 463--482.

\bibitem[Ukk92]{Ukkonen}
E.~Ukkonen, \emph{Approximate string-matching with q-grams and maximal
  matches}, Theoretical Computer Science \textbf{92} (1992), 191--211.

\bibitem[Wel93]{Welsh93}
D.J.A. Welsh, \emph{Complexity: Knots, colourings, and counting}, London
  Mathematical Society Lecture Notes Series, no. 186, Cambridge University
  Press, Cambridge, England, 1993.

\bibitem[Wil90]{Wilf}
Herbert~S. Wilf, \emph{generatingfunctionology}, 2nd ed., Academic Press, San
  Diego, CA, 1990.

\end{thebibliography}
\end{document}